\documentclass[final]{ws-ccm}
\usepackage{mathtools}
\usepackage{mathrsfs}

\newcommand{\LL}{\mathop{\hbox{\vrule height 7pt width .5pt depth 0pt
\vrule height .5pt width 6pt depth 0pt}}\nolimits}

\newcommand{\res}{\mathop{\hbox{\vrule height 7pt width .5pt depth 0pt
\vrule height .5pt width 6pt depth 0pt}}\nolimits}

\newcommand{\Om}{\Omega}
\newcommand{\eps}{\varepsilon}
\newcommand{\ep}{\varepsilon}
\newcommand{\DD}{{\cal D}}
\renewcommand{\AA}{{\mathcal A}}
\newcommand{\BB}{{\mathcal B}}
\newcommand{\GG}{{\mathcal G}}

\newcommand{\om}{\omega}

\newcommand{\curl}{{\rm curl\, }}
\newcommand{\grad}{\nabla}

\newcommand{\dist}{{\rm dist\, }}

\newcommand{\RR}{{\mathbb R}}
\newcommand{\CC}{{\mathbb C}}
\newcommand{\ZZ}{{\mathbb Z}}
\newcommand{\NN}{{\mathbb N}}

\newcommand\lep{|\ln \ep|}
\renewcommand{\theequation}{\thesection .\arabic{equation}}
\newcommand{\be}{\begin{equation}}
\newcommand{\ee}{\end{equation}}
\newcommand{\bea}{\begin{eqnarray}}
\newcommand{\eea}{\end{eqnarray}}
\newcommand{\beann}{\begin{eqnarray*}}
\newcommand{\eeann}{\end{eqnarray*}}

\begin{document}

\markboth{S. Alama,  L. Bronsard \&  V. Millot} {$\Gamma$-convergence of 2D Ginzburg-Landau functionals}

\title{$\Gamma$--CONVERGENCE OF 2D GINZBURG-LANDAU FUNCTIONALS\\ 
WITH VORTEX CONCENTRATION ALONG CURVES}

\author{STAN ALAMA}
\address{Department of Mathematics and Statistics\\
McMaster University\\
Hamilton, Ontario, Canada L8S 4K1\\
\emph{\tt{alama@mcmaster.ca}}}

\author{LIA BRONSARD}
\address{Department of Mathematics and Statistics\\
McMaster University\\
Hamilton, Ontario, Canada L8S 4K1\\
\emph{\tt{bronsard@mcmaster.ca}}}

\author{VINCENT MILLOT}
\address{Universit\'e Paris Diderot - Paris 7\\
CNRS, UMR 7598 Laboratoire J.L. Lions\\
F-75005 Paris, France\\
\emph{\tt{millot@math.jussieu.fr}}}

\maketitle

\begin{abstract}
{\bf Abstract.}  
We study the variational convergence 
of a family of two-dimensional Ginzburg-Landau functionals arising in the study of superfluidity or thin-film superconductivity, as the Ginzburg-Landau parameter $\eps$ tends to $0$. In this regime 
and for large enough applied rotations (for superfluids) or magnetic fields (for superconductors), 
the minimizers acquire quantized point singularities (vortices).  We focus on situations in which an unbounded number of vortices accumulate along a prescribed Jordan curve or a simple 
arc in the domain.  This is known to occur in a circular annulus under uniform rotation, or in a simply connected domain with an appropriately chosen rotational vector field.   We prove that, suitably normalized, the energy functionals $\Gamma$-converge to a classical energy from potential theory.  Applied to global minimizers, our results describe the limiting distribution of vortices along the curve in terms of Green equilibrium measures.
\end{abstract}



\section{Introduction}                        

The Ginzburg-Landau theories have had an enormous influence on both physics and mathematics.  
Physicists employ Ginzburg-Landau models in modeling superconductivity, superfluidity, and, more recently, for rotating Bose-Einstein condensates (BECs), 
all systems which present quantized defects commonly known as vortices.  In mathematics, starting with the work by Bethuel, Brezis \& H\'elein \cite{BBH2}, 
many powerful methods have been developed to study the physical London limit, {\it i.e.}, as the Ginzburg-Landau parameter $\varepsilon$ tends to $0$.  
This limit corresponds to the Thomas-Fermi regime in BEC, and to an analogous regime in superfluids where the characteristic length scale $\eps$ is very small.  
In a two-dimensional setting,  vortices are essentially characterized as isolated zeroes of the order parameter carrying a winding number, and in the London limit
as point defects where energy concentration occurs.  
The question of whether energy minimizers develop vortices, where they appear in the domain, and how many there should be (for given boundary conditions, 
constant applied fields or angular velocities) has been analyzed in many contexts and parameter regimes.
\vskip5pt

In this paper, we will focus on the following Ginzburg-Landau energy, arising for instance in the physical context of a rotating superfluid. Considering 
a bounded simply connected domain $\DD\subset\RR^2$,  a smooth vector field $V:\RR^2\to\RR^2$, $\Omega>0$ and $\eps>0$,   we define the functional
$$  u\in H^1(\DD;\CC)\mapsto   {\mathcal F}_\eps(u):= \int_\DD \left\{  \frac12|\nabla u|^2 
               + \frac{1}{4\eps^2}(1-|u|^2)^2-\Omega\, V(x)\cdot j(u) \right\} dx\,.  $$
Identifying $\RR^2$ with the complex plane $\CC$, we denote by
$$j(u):=u\wedge\nabla u\in L^1(\DD;\RR^2)\,,$$ 
the {\it pre-Jacobian} of $u$. The $L^1$-vector field $j(u)$ may also be written as $j(u)=(iu,\nabla u)$, 
where $(\cdot,\cdot)$ is the standard inner product of two complex numbers, viewed as vectors in~$\RR^2$.
\vskip5pt

In the case of uniform rotation, that is $V(x)=x^\perp=(-x_2,x_1)$ and with $\DD$ a disk, 
Serfaty  \cite{S-Super}   studied minimizers of a closely related functional (see Remark~\ref{relfunct}) to determine the critical value $\Omega_1=\Omega_1(\varepsilon)$ of the angular speed $\Omega$ 
at which vortices first appear (see also \cite{IM1,IM2} for BECs).  She finds that
minimizers acquire vorticity at $\Omega_1= k(\DD)|\ln\eps| + O(\ln|\ln \eps|)$ for an explicitly determined  constant $k(\DD)$. 
In a series of papers, culminating with the publication of the research monograph \cite{SSbouquin}, Sandier \& Serfaty developed powerful tools to study vortices in 
Ginzburg-Landau models.  Although they primarily work with the full Ginzburg-Landau model with magnetic field, the methods apply as well to the functional ${\mathcal F}_\eps$ above.  
In particular, their results apply to the near-critical regime in {\em simply connected} domains.  In our setting, their results show that for any simply connected domain $\DD$, the first order 
expansion of the critical value $\Omega_1$ for vortex existence in minimizing configurations is also given by $k(\DD)|\ln\eps|$ for some constant  $k(\DD)$.  
Moreover the locus of concentration of vortices for $\Omega=\Omega_1 + o(|\ln\eps|)$ is given by the set of maxima of $|\zeta|$, with $\zeta$ the solution of the following boundary-value problem:
\begin{equation}\label{detzeta}
\begin{cases}
-\Delta \zeta = \curl V & \text{in $\DD$}\,, \\
\zeta=0 & \text{on $\partial\DD$}\,.
\end{cases}   
\end{equation}
The constant $k(\DD)$ is then determined by
$$k(\DD) =\frac{1}{2|\zeta|_{\max}} \,,$$
where $|\zeta|_{\max}$ denotes the maximum value of  $|\zeta|$. If, for instance, $V$ is real-analytic and $\curl V$ is nonnegative, then so is the solution $\zeta$, and the maximum is  generically 
attained at a finite number of points in $\DD$ (see {\it e.g.} \cite{BJS}). 
In this situation, if  $\Omega=\Omega_1 + o(|\ln\eps|)$, minimizers exhibit {\em vortex concentration at isolated points}, and the number of vortices remains uniformly bounded whenever 
$\Omega-\Omega_1$ is of order $O(\ln|\ln\varepsilon|)$, see \cite{SSbouquin, S-Super,IM1,IM2}. 
\vskip5pt

The case of a multiply connected domain provides a slightly different qualitative picture.  
In a work on rotating Bose-Einstein condensates, Aftalion, Alama \& Bronsard~\cite{AAB} considered a similar functional 
in a domain given by a circular annulus $\mathcal{A}$ (centered at the origin) and again with uniform rotation $V(x)=x^\perp$ (see Remark~\ref{relfunct}). 
Unlike the simply connected case, minimizers in the annulus may have vorticity without vortices, as the hole acquires positive winding at bounded rotation $\Omega$.  
Then point vortices are nucleated inside the interior of $\mathcal{A}$ at a critical value $\Omega_1$, again of leading order $|\ln\eps|$.  
Solving equation \eqref{detzeta} in the annulus $\mathcal{A}$, one finds out that the set of maxima of the function 
$\zeta$ is given by a circle inside $\mathcal{A}$ (see Example~\ref{exhole}). Hence one can expect that, rather than accumulating at isolated points, vortices concentrate along  this circle 
in the limit $\eps\to 0$.  
The main feature proved in \cite{AAB} is that if $\Omega\sim\Omega_1+ O(\ln|\ln\eps|)$, then vortices are indeed essentially supported by 
a circle $\Sigma$ and that the total degree of these vortices is of order $\ln|\ln\eps|$. In other words, in the limit $\varepsilon\to 0$,  infinitely many vortices concentrate on $\Sigma$, 
a phenomenon that we call {\em vortex concentration along a curve}. However the question of  the distribution of the limiting vorticity around the circle was left open. 
Subsequent results of Alama \& Bronsard \cite{AB1,AB2} extend the result of \cite{AAB} to multiply connected domains and to the full Ginzburg-Landau model with magnetic field 
and pinning potential. In contrast with the previous case, concentration on curves might not be a generic phenomenon for the Ginzburg-Landau model with magnetic field in a general multiply connected domain. 
Indeed, in this setting the vector field $V$ represents the electromagnetic potential and it is an unknown of the problem. The results in \cite{AB1,AB2} show that, for near-critical external applied fields, 
the concentration set of vortices is also given by  the set of maxima of a certain potential  related to $V$. This set may contain finitely many points and/or closed loops. Assuming that it  
contains closed loops,  they prove  that {\em vortex concentration along a curve} occurs, but the determination of the limiting vorticity was again left open. 
\vskip5pt

To effectively separate the question of the nature of the concentration set from the question of localizing vortices, we instead start with a 
{\em simply connected} domain $\DD$, and we prescribe the function $\zeta$ with $\zeta\ge 0$ in $\DD$ and  $\zeta|_{\partial\DD}=0$, in such 
a way that $\zeta$ is maximized on a {\em prescribed curve} $\Sigma\subset\!\subset \DD$.  Then, we define 
$$  V(x):=-\nabla^\perp\zeta(x)=\left(\frac{\partial\zeta}{\partial x_2}, -\frac{\partial\zeta}{\partial x_1}\right)$$ 
as our vector field.
We will prove that vortices will be forced to accumulate on $\Sigma$ as $\eps$ tends to $0$.
The curve $\Sigma$ can be either a smooth Jordan curve or a smooth embedded simple arc, compactly contained in $\DD$.  In this setting, we shall resolve the problem of distribution of vortices along curves, 
both for minimizers and in the more general setting of $\Gamma$-convergence. In the last section we will show that in a multiply connected domain, and for more general vector fields $V$,  
the problem does not differ too much in nature, and that a similar analysis can be performed. 
\vskip5pt

To state our main result we must give more specific hypotheses on $\zeta$ and the angular speed $\Omega$.  We  assume that $\zeta$ satisfies the following assumptions:
\begin{itemize}
\item[(H1)]\label{H1} $\zeta\in C^{0,1}_0(\DD)$, $\zeta\ge 0$ in 
$\overline{\DD}$,  and $\zeta_{\max}:=\max_{x\in\overline{\DD}}\, \zeta(x)>0$;
\vskip5pt

\item[(H2)]\label{Sig} $\Sigma:=\{x\in\DD: \ \zeta(x)=\zeta_{\max}\}\subset\!\subset \DD$ is
a Jordan curve or a simple embedded arc of class $C^2$.
 \end{itemize}
We further assume that $\Omega=\Omega(\eps)$ is near to the critical value needed for the presence of vortices.  More precisely, 
\begin{equation}\label{omega}
\Om_\varepsilon:=\frac{\lep}{2\zeta_{\rm max}} + \om(\varepsilon)\,,
\end{equation}
for some function $\omega:(0,+\infty)\to(0,+\infty)$  satisfiying $\omega(\varepsilon)\to+\infty$ with $|\ln\varepsilon|^{-1}\omega(\varepsilon)\to 0$ as $\varepsilon\to 0$.
\vskip5pt

For $u\in H^1(\DD;\CC)$ we consider the rescaled functional
$$  F_\eps(u) := \frac{1}{\omega^2(\eps)}\int_\DD
   \left\{  \frac12 |\nabla u|^2 + \frac{1}{4\eps^2} (1-|u|^2)^2+ \Omega_\eps \nabla^\perp\zeta\cdot j(u) \right\} dx\,,$$
and for a nonnegative Radon measure $\mu$ on $\mathcal{D}$, we define 
$$  I(\mu):=\frac{1}{2}\iint_{\DD\times \DD}
G(x,y) \, d\mu(x)\, d\mu(y) \,, $$
where the function $G$ denotes the Dirichlet Green's function of the domain $\mathcal{D}$, {\it i.e.}, for every $y\in\mathcal{D}$, $G(\cdot,y)$ is the solution of 
\begin{equation}\label{defGreen}
\begin{cases}
-\Delta G(\cdot,y)=\delta_y &\text{in $\mathscr{D}'(\mathcal{D})$}\,,\\
G(\cdot,y) =0 & \text{on $\partial\mathcal{D}$}\,.
\end{cases}
\end{equation}
\vskip5pt

Our main result  deals with the $\Gamma$-convergence of the family of functionals $\{F_\varepsilon\}_{\varepsilon>0}$ as $\varepsilon\to 0$, and it is stated (as usual)
in terms of the vorticity distribution given by the weak Jacobian, that is half the distributional curl of the pre-Jacobian (see  {\it e.g.} \cite{SSbouquin}).

\begin{theorem}\label{mainthm}
Assume that \emph{(H1)}, \emph{(H2)} and \eqref{omega} hold. 
Let $\varepsilon_n\to 0^+$ be an arbitrary sequence. Then,
\begin{itemize}
\item[(i)]  for any $\{u_n\}_{n\in\NN}\subset H^1(\mathcal{D};\CC)$ satisfying $\sup_n F_{\varepsilon_n}(u_n)<+\infty$, there exist a subsequence (not relabelled) and a 
nonnegative Radon measure $\mu$ in $H^{-1}(\mathcal{D})$ supported by $\Sigma$ such that 
\begin{equation}\label{convjac}
\frac{1}{\omega(\varepsilon_n)}\, \curl\,j(u_n)\mathop{\longrightarrow}\limits_{n\to+\infty} \mu \quad\text{strongly in $(C^{0,1}_0(\mathcal{D}))^*$}\,;
\end{equation}
\item[(ii)] for any $\{u_n\}_{n\in\NN}\subset H^1(\mathcal{D};\CC)$ such that \eqref{convjac} holds for some nonnegative Radon measure $\mu$ in 
$H^{-1}(\mathcal{D})$ supported by $\Sigma$, we have 
$$\liminf_{n\to+\infty}\,F_{\varepsilon_n}(u_n)\geq I(\mu)-\zeta_{\rm max}\,\mu(\mathcal{D})\,;$$
\item[(iii)] for any  nonnegative Radon measure $\mu$ in $H^{-1}(\mathcal{D})$ supported by $\Sigma$, there exists a sequence $\{u_n\}_{n\in\NN}\subset H^1(\mathcal{D};\CC)$ such that 
\eqref{convjac} holds and 
$$\lim_{n\to+\infty} F_{\varepsilon_n}(u_n)= I(\mu)-\zeta_{\rm max}\,\mu(\mathcal{D})\,.$$
\end{itemize}
\end{theorem}

As it is well known, the $\Gamma$-convergence theory is well suited to study  asymptotics in minimization problems  (see {\it e.g.} \cite{DalMaso}). 
In this context, we shall derive from Theorem~\ref{mainthm} the following convergence result for the vorticity of global minimizers, and hence solving 
the problem  on the limiting distribution of vortices along $\Sigma$, see Remark \ref{equimeas} below. 

\begin{corollary}\label{asymptmin}
Assume that \emph{(H1)}, \emph{(H2)} and \eqref{omega} hold. 
Let $\varepsilon_n\to 0^+$ be an arbitrary sequence. For every integer $n\in\mathbb{N}$, let $u_n\in H^1(\mathcal{D};\CC)$ be a minimizer of $F_{\varepsilon_n}$. 
Then, 
$$\frac{1}{\omega(\varepsilon_n)}\, \curl\,j(u_n)\mathop{\longrightarrow}\limits_{n\to+\infty} \,\frac{\zeta_{\rm max}}{2 I_*}\,\mu_* \quad\text{strongly in $(C^{0,1}_0(\mathcal{D}))^*$}\,,$$
where $\mu_*$ is the unique minimizer of $I$ over all  probability measures supported on $\Sigma$, and $I_*:=I(\mu_*)$.
\end{corollary}

\begin{remark}\label{bestn}
As a direct application of the results in \cite{SSbouquin} (see also \cite{JS}), we shall see in Section 2  that for  configurations $\{u_\varepsilon\}$ with $F_\varepsilon$-energy uniformly bounded from above, 
the vorticity distribution 
$\curl j(u_\varepsilon)$ can be  approximated (with respect to the $(C^{0,1}_0(\mathcal{D}))^*$--topology) by a measure of the form $2\pi\sum_{i\in I_\varepsilon}d_i\delta_{a_i}$ for some finite set of points 
$\{a_i\}_{i\in I_\varepsilon}\subset \DD$ and integers $\{d_i\}_{i\in I_\varepsilon}\subset \mathbb{Z}$.  
In other words, each point $a_i$ can be viewed as an approximate vortex with winding number $d_i$. Thus the integer $D_\varepsilon =\sum_{i\in I_\varepsilon} |d_i|$ may be refered to as to 
{\it approximate total  vorticity} of the configuration $u_\varepsilon$. 
It is commonly known that approximate vortices carry a kinetic energy essentially greater than or equal to 
$\pi D_\varepsilon|\ln\varepsilon|$ (see Section~2 for more details). With such an estimate in hand, and 
using the arguments of Section~3, we actually obtain a more refined lower bound for the energy than the one given by Theorem~\ref{mainthm}, claim {\it (ii)}. More precisely, one has
$$\liminf_{\eps\to 0} \frac{1}{\omega^2(\eps)} \left( \int_\DD \frac12 |\nabla u_\eps|^2 + \frac{1}{4\eps^2}(1-|u_\eps|^2)^2\,dx - \pi D_\eps \lep \right) \ge I(\mu)\,,$$
and
$$\liminf_{\eps\to 0} \frac{1}{\omega^2(\eps)} \left( \Omega_\eps\int_\DD \nabla^\perp\zeta\cdot j(u_\eps)\,dx + \pi D_\eps\lep\right) \ge -\zeta_{\rm max}\, \mu(\DD)\,.
$$
As a consequence, if $\{u_\eps\}$ is any recovery sequence (in the sense of {\it (iii)} of Theorem~\ref{mainthm}), the $\liminf$'s above become limits, and equality holds in each case.  In analogy with \cite{BBH2},  we may then say that $I(\mu)$ plays the role of  renormalized energy.
\end{remark}

\begin{remark}[Minimizers]
From Corollary \ref{asymptmin} and Remark~\ref{bestn}, we deduce that if $u_\varepsilon$ is energy minimizing, then 
$$   D_\varepsilon = \frac{\zeta_{\rm max}}{4\pi I_*}\, \omega(\varepsilon)+o(\omega(\varepsilon)) \quad\text{as $\varepsilon\to 0$}\,, $$  
and from Theorem~\ref{mainthm}, the minimal value of the energy  expands as 
$$   \min_{H^1(\DD;\CC)}\, \omega^2(\varepsilon)F_{\varepsilon}=  
-\frac{\zeta^2_{\rm max}}{4I_*}\,\omega^2(\eps)+ o(\omega^2(\varepsilon))\,. $$
\end{remark}

\begin{remark}[Equilibrium measures] \label{equimeas}
The value $I(\mu)$ gives the electrostatic energy of a positive charge distribution $\mu$ on the set $\Sigma\subset\!\subset\DD$.  
The minimizer $\mu_*$ of $I$ over all probability measures on $\Sigma$ is called the {\em Green equilibrium measure} in $\DD$ associated to the set 
$\Sigma$, and gives the equilibrium charge distribution of a charged conductor inside of a neutral conducting shell, represented by $\partial\DD$. 
The  value $1/I_*$ is refered to as to the {\em capacity of the condenser} $(\Sigma,\partial \DD)$.  
The  interested reader can find in \cite{ST} many results on the existence and general (regularity) properties of the equilibrium measures as well as some examples. 
For instance, if $\DD$ is a disc and $\Sigma$ is a concentric circle, 
then the equilibrium measure $\mu_*$ is the normalized arclength measure on $\Sigma$, see \cite[Example~II.5.13]{ST}, and thus vortices are asymptotically equidistributed along $\Sigma$ as $\eps\to 0$.  
However for an arbitrary curve $\Sigma$, the distribution is of course non-uniform in general. In case where $\Sigma$ is an embedded arc, it is  even singular at the endpoints, see \cite[Example~II.5.14]{ST}. 
\end{remark}

\begin{remark}\label{diffsigma} 
In the present results the structure and regularity assumptions on the set $\Sigma$ given in (H2) are  
mainly motivated by the physical context of \cite{AAB,AB1,AB2}. However it will be clear  that (H2) can be 
relaxed into weaker statments. More precisely, the proof of Theorem~\ref{mainthm} relies on (H2) only for the $\Gamma-\limsup$ inequality, {\it i.e.}, claim {\it (iii)}. 
The construction of the recovery sequence (see Section~\ref{secupbd}) could be applied with minor modifications if the set $\Sigma$ is for instance a finite union of piecewise $C^2$ arcs/Jordan curves. 
Actually $\Sigma$ could even have a more general structure such as a non-empty interior. In this later case we assume that $\partial\Sigma$ is made by 
finitely many arcs and Jordan curves of class $C^2$. Then, given a nonnegative Radon measure $\mu\in H^1(\DD)$ supported by $\Sigma$, one can construct a recovery sequence for $\mu$ 
applying the approximation techniques of Section~\ref{secupbd} to  $\mu\res\partial\Sigma$ and a more standard regularization procedure for $\mu\res{\rm int}(\Sigma)$ as in \cite{SSbouquin}. 
\end{remark}

In Section~\ref{annulus} we will show how to apply the method to a multiply connected domain and a more general vector field $V$.  
For simplicity, we shall consider only domains $\mathcal{A}$ which are topological annuli, {\it i.e.},  
$\mathcal{A}=\DD\setminus\overline\BB$, where $\DD$, $\mathcal{B}$ are simply connected and $\BB\subset\!\subset\DD$.  
For multiply connected domains and/or a general field $V$, 
there is an extra step involved in the analysis. Indeed, the leading order term in the minimal energy (of order $|\ln\eps|^2$) is due 
to the curl-free part of $V$ which induces a diverging phase in any minimal configuration,  
and to the vorticity in the hole $\BB$ which acts as a sort of giant vortex (by analogy with \cite{AAB}).  In other words, to get information on the internal vorticity,  
one has to perform a second order $\Gamma$-convergence analysis.  
As in \cite{AB1} we first describe minimizing vortexless configurations ({\it i.e.}, energy minimizers over $\mathbb{S}^1$-valued maps) which nonetheless have vorticity around the hole. 
Then we show that for arbitrary configurations, the energy due to the curl-free part of $V$ and the hole decouples nearly exactly (see Proposition~\ref{energdecomp}). 
After separating out this contribution to the energy, 
the residual energy functional resembles $F_\eps$ above, and the $\Gamma$-convergence analysis can be done in much the same way, with some care taken to control the residual vorticity around the hole. 
The result is stated in Theorem~\ref{thmhole}.  The $\Gamma$-limit again involves a Green energy for measures supported by a prescribed set.

\begin{remark}[Related funtionals]\label{relfunct}
As mentioned earlier, the methods here may also be applied to other Ginzburg-Landau functionals which have the same structure as ${\mathcal F}_\eps$.
A simple variant is
$$\tilde {\mathcal F}_\eps (u) = \int_\DD
   \left\{  \frac12 |(\nabla - i \Omega_\eps V)u|^2  + \frac{1}{4\eps^2} (1-|u|^2)^2\right\} \,dx\,.$$
This energy has been studied in more than one context.  Serfaty considers in \cite{S-Super}  the minimization of this energy for rotating superfluids, with $V=x^\perp$ 
and under a Dirichlet boundary condition $u|_{\partial\DD}=0$.  The minimization of the same energy under natural (Neumann) boundary conditions also arises in a 
simplified model of thin-film superconductors introduced by Chapman, Du \& Gunzburger \cite{CDG} (see also Alama, Bronsard \& Galv\~ao-Sousa \cite{ABGS}). 
In this setting, $\Omega_\eps V(x)=A_\eps(x)$ represents the magnetic vector potential of the externally applied magnetic field $h_{ex}=\Omega_\eps\curl V$. 
With $\Omega_\eps$ satisfying \eqref{omega}, and a sequence $\{u_\varepsilon\}$ in the energy regime of Theorem~\ref{mainthm} (which holds for minimizers under a 
homogeneous Neumann boundary condition), 
the two energies agree very closely, 
$$  \tilde {\mathcal F}_\eps(u_\varepsilon) ={\mathcal F}_\eps(u_\varepsilon)+\frac{\Omega^2_\varepsilon}{2}\int_\DD|V|^2\,dx+o(1)\quad\text{as $\varepsilon\to 0$}\,.$$
In the case of the homogeneous Dirichlet boundary condition some care must be taken since 
 a singular boundary layer arises near $\partial \DD$ as $\varepsilon\to 0$ (see \cite{S-Super}), but otherwise the result of Theorem~\ref{mainthm} should remain essentially the same.

For annular domains, another functional which resembles ${\mathcal F}_\eps$ has been used in the modeling of rotating BECs in certain anharmonic traps (see \cite{AAB}),
$$  \GG_\varepsilon(u) = \int_\AA
   \left\{  \frac12|\nabla u|^2
              + \frac{1}{4\eps^2} (a(x)-|u|^2)^2-\Omega_\varepsilon\, V(x)\cdot j(u)   \right\} dx\,.
$$
Here $V(x)=x^\perp$ is the velocity field of uniform rotation, and the function $a(x)$, positive in $\AA$,  gives the trapping potential which contains the condensate. 
It is shown in \cite{AAB} that $\mathcal{G}_\varepsilon$-minimizers develop as $\varepsilon\to 0$, infinitely many vortices concentrating along a circle for sufficently high rotation $\Omega_\varepsilon$. 
We believe that similar results to the ones in Section~\ref{annulus} should hold, but the method is not directly applicable here. Indeed, our analysis is 
based on ``global energy estimates" of \cite{SSbouquin} (see \eqref{balls3} in Proposition \ref{consvortexball}) and the presence of the inhomogeneity $a(x)$ requires local estimates.  
Moreover the analysis of vortices for the energy $\GG_\varepsilon$ is complicated 
by the fact that $a(x)$ vanishes on the boundary, and some delicate estimates are required so as not to lose too much information near the boundary.

One may also consider the full Ginzburg-Landau model of superconductivity for complex order parameter $u: \DD\to\CC$ and 
magnetic vector potential $A:\DD\to\RR^2$.  For the Ginzburg-Landau model in a simply connected domain 
with constant applied magnetic field, Sandier \& Serfaty \cite[Theorem 9.1]{SSbouquin} have proven a $\Gamma$-convergence theorem of the form of 
Theorem~\ref{mainthm} for applied fields of the form $h_{ex}=H_{C1} + \omega(\eps)$, where in that case it is appropriate to take 
$\ln|\ln\eps|\lesssim \omega(\eps) \ll|\ln\eps|.$  This problem exhibits {\it vortex concentration at points}, and the limiting energy is obtained 
by rescaling around  the points of concentration.  The rescaled vorticity measures of minimizers converge to an equilibrium measure associated to a different problem in 
potential theory, a ``Gauss variation'' problem whereby the charges are to be optimally placed in $\RR^2$ subject to an applied electric field (see \cite{ST}). 
Concentration on curves is possible in multiply connected regions (see \cite{AB1,AB2}). The methods described here should apply in this situationl once adapted to the magnetic setting, 
although the energy obtained as a $\Gamma$-limit should be for the Helmholtz (and not the Laplace) Green's function, $-\Delta_x H(x,y) + H(x,y)=\delta_y(x)$.  
However, for the full Ginzburg-Landau functional with magnetic field in general domains, it is an interesting open problem in PDE to determine which (if any) non-symmetric multiply connected 
domains exhibit vortex concentration on curves. 

We finally mention a recent paper by Kashmar \cite{Kash} exhibiting concentration on a circle for a magnetic Ginzburg-Landau functional in the disc with 
an inhomogeneity $a(x)$ as above described by a radial step function (modelling for instance a superconducting body made by two different species).  Here again we believe 
that  similar results should hold but the method does not directly apply due to the inhomogeneity of the Ginzburg-Landau energy density. 
\end{remark}

The plan of the paper is as follows.  In Section 2 we prove assertion {\it (i)} of Theorem~\ref{mainthm}.  
Section~3 tackles part {\it (ii)} of the theorem. The upper bound of statement {\it (iii)} is derived in Section~4, completing the proof of Theorem~\ref{mainthm}.  
The proof of Corollary~\ref{asymptmin} is presented at the end of Section~4. 
Section~\ref{annulus} sketches how the preceeding arguments must be modified to treat annular domains.  
\vskip10pt

\noindent{\bf Notations.}
For any open set $B\subset\DD$ and any admissible map $u$, we denote by 
$$E_{\varepsilon}(u,B):= \int_B \left\{  \frac12|\nabla u|^2 
               + \frac{1}{4\eps^2}(1-|u|^2)^2 \right\} \,dx$$
the so-called Ginzburg-Landau energy of $u$ in $B$. Given $\varepsilon>0$ we define 
$$  \DD_\eps:=\{x\in \DD: \  \dist(x,\partial\DD)>\eps\}\,,  $$
and for a sequence $\varepsilon_n\to0^+$, we shall write $\omega_n:=\omega(\varepsilon_n)$ and $\Omega_n:=\Omega_{\varepsilon_n}$.

\section{Compactness of normalized weak Jacobians}\label{degbounded}                               

This section is devoted to the proof of claim {\it (i)} in Theorem \ref{mainthm}. 
The key ingredient to prove compactness of normalized weak Jacobians is the so-called  {\it ``vortex balls construction"} taken from \cite[Theorem~4.1]{SSbouquin}. 

\begin{proposition}\label{consvortexball}
There exists a constant $\varepsilon_0>0$ such that for any $0<\varepsilon<\varepsilon_0$ and any $u\in H^1( \mathcal{D};\mathbb{C})\cap C^1( \mathcal{D})$ 
satisfying $E_\varepsilon(u,\mathcal{D})\leq \sqrt{\varepsilon}$, the following holds. For any $C_0\varepsilon^{1/4}<r<1$ there exists a finite collection of disjoint closed balls 
$\mathcal{B}_r=\{\overline B(a_i,\rho_i)\}_{i\in I_r}$ such that, writing $B_i:=\overline B(a_i,\rho_i)$,  
\begin{itemize}
\item[(i)] $\displaystyle r=\sum_{i\in I_r}\rho_i\,$;
\vskip5pt
\item[(ii)] $\displaystyle \big\{|1-|u||\geq \varepsilon^{1/8}\big\}\cap \mathcal{D}_\varepsilon \subset\mathcal{D}_\varepsilon\cap \bigcup_{i\in I_r}B_i=:V\,$;\
\vskip5pt
\item[(iii)] setting $d_i:={\rm deg}(u,\partial B_i)$ if $B_i\subset\mathcal{D}_\varepsilon$, and $d_i:=0$ otherwise,  
\begin{equation}\label{balls3}
E_\varepsilon(u,V)\geq \pi D_r\bigg(\ln\bigg(\frac{r}{\varepsilon D_r}\bigg)-C_1\bigg) \,,
\end{equation}
where $\displaystyle D_r:=\sum_{i\in I_r}|d_i|\,$ is assumed to be positive;
\vskip5pt
\item[(iv)] $\displaystyle D_r\leq C_2|\ln\varepsilon|^{-1} E_\varepsilon(u,\mathcal{D})\,$;
\end{itemize}
and $C_0$, $C_1$, $C_2$ are universal  constants. Moreover, if $C_0\varepsilon^{1/4}<r_1<r_2<1$ and $\mathcal{B}_1$, $\mathcal{B}_2$ are the corresponding families of balls, 
then every ball of $\mathcal{B}_1$ is included in a ball of $\mathcal{B}_2$.
\end{proposition}

In the remainder of this section, we consider an arbitrary sequence $\varepsilon_n\to0^+$. In the following lemma, we  prove an upper bound on the Ginzburg-Landau energy 
for a sequence having an $F_{\varepsilon_n}$-energy uniformly bounded from above. It will allow us to apply the previous proposition to such a sequence. 

\begin{lemma}\label{lem1} Assume that \emph{(H1)} and \eqref{omega} hold.  
Let $\{u_n\}_{n\in\NN}$ be a sequence in $H^1(\DD;\mathbb{C})$  such that $\sup_n F_{\varepsilon_n}(u_n)<+\infty$.  Then there exists a constant $C$ independent of $n$ such 
that $\|u_n\|_{L^4(\mathcal{D})}\leq C$ and $E_{\varepsilon_n}(u_n,\mathcal{D})\le C |\ln\varepsilon_n|^2$. 
\end{lemma}

\begin{proof}
Observe that 
\begin{multline}\label{contrGL}
E_{\varepsilon_n}(u_n,\mathcal{D})= \omega^2_nF_{\varepsilon_n}(u_n) - \Om_{n}\int_\DD \nabla^\perp\zeta\cdot j(u_n)\, dx
\leq  \Om_{n}\|\nabla\zeta\|_\infty\int_\DD|u_n||\nabla u_n|dx +O(\omega_n^2) \\
\leq  \frac{1}{4}\int_{\DD}|\nabla u_n|^2dx+\Omega_{n}^2\|\nabla\zeta\|_\infty^2\int_{\DD}|u_n|^2 dx +O(\omega_n^2)\,.
 \end{multline}
In particular, 
$$\int_{\DD}\big(|u_n|^4-2(1+2\varepsilon_n^2\Omega_{n}^2\|\nabla\zeta\|_\infty^2)|u_n|^2+1\big)dx \leq O(\varepsilon^2_n\omega_n^2)\,,$$
so that $\|u_n\|_{L^4(\DD)}\leq C$ for a constant $C$ independent of $n$. Inserting this estimate in \eqref{contrGL}, the announced result follows easily. 
\end{proof}

The first step in proving compactness of the normalized Jacobians is to show that the {\it approximate total vorticity} is bounded by the excess rotation $\om_n$. 
We emphasize that here $\Sigma$ could be any compact subset of $\DD$.  

\begin{proposition}\label{degreebound} Assume that \emph{(H1)} and \eqref{omega} hold.    
Let $\{u_n\}_{n\in\NN}$ be a sequence in $H^1(\DD;\mathbb{C})\cap C^1(\mathcal{D})$  such that $\sup_n F_{\varepsilon_n}(u_n)<+\infty$. Then there exist $r_n\to 0^+$ 
and a sequence of families of balls $\mathcal{B}_{r_n}=\{\overline B(a_i^n,\rho_{i,n})\}_{i\in I_{r_n}}$ as in Proposition \ref{consvortexball} 
such that $D_{r_n}\leq C\omega_n$ for some constant $C>0$ independent of $n$. 
\end{proposition}

\begin{proof}  
 Let  $r_n:=|\ln\varepsilon_n|^{-4}$. In view of Lemma \ref{lem1} we can apply Proposition \ref{consvortexball} to $u_n$ with $r=r_n$ and $n$ large enough. 
 For each such $n$ we denote by $\mathcal{B}_{r_n}=\{\overline B(a_i^n,\rho_{i,n})\}_{i\in I_{r_n}}$ the corresponding family of balls. For convenience we write 
$I_n:=I_{r_n}$, $B_i^n:=\overline B(a_i^n,\rho_{i,n})$, $d_{i,n}:={\rm deg}(u,\partial B^n_i)$ if $B_i^n\subset \mathcal{D}_{\varepsilon_n}$, $d_{i,n}=0$ otherwise, and $D_n:=D_{r_n}$. 
\vskip5pt
 
\noindent{\it Step 1.} We first claim that 
\begin{equation}\label{lem2}
\Omega_{n}\int_\DD \nabla^\perp\zeta\cdot j(u_n)\, dx
    = -2\pi\Omega_{n}\sum_{i\in I_{n}} d_{i,n}\, \zeta(a^n_i) + o(1)\quad\text{as $n\to+\infty$}\,.
\end{equation}
Indeed, applying \cite[Theorem~6.1]{SSbouquin} with the balls $\{B_i^n\}_{i\in I_n}$ 
we derive the estimate
\begin{equation}\label{jacobian}  
\bigg\|  \curl j(u_n) - 2\pi \sum_{i\in I_n} d_{i,n}\, \delta_{a^n_i} \bigg\|_{(C_0^{0,1}(\mathcal{D}))^*}
   \leq C r_n \big(1+E_{\varepsilon_n}(u_n,\mathcal{D})\big) \leq C|\ln \varepsilon_n|^{-2}\,,  
\end{equation}
thanks to Lemma \ref{lem1}. Then \eqref{lem2} follows since $\Omega_n=O(|\ln\varepsilon_n|)$ and $\zeta\in C_0^{0,1}(\DD)$.
\vskip5pt

\noindent{\it Step 2.} Without loss of generality we may assume that $\omega_n=O(D_n)$, otherwise there is nothing to prove. 
In view of claim {\it (iv)} in Proposition \ref{consvortexball} and Lemma \ref{lem1}, we have $D_n\leq O(|\ln\varepsilon_n|)$. Next, combining the 
lower bound (\ref{balls3}) with \eqref{lem2}, we infer that 
\begin{multline}\label{lb1}
O(\omega^2_n)\geq  \omega^2_n F_{\varepsilon_n}(u_n)
\geq   \pi D_{n} \big(|\ln\varepsilon_n|-C\ln|\ln\varepsilon_n|\big)
     - 2\pi\Omega_{n}\sum_{i\in I_{n}} d_{i,n}\zeta(a^n_i) \,+\\ 
+  \int_{\DD\setminus\cup_{i\in I_n} B^n_i} |\nabla u_n|^2\, dx +o(1)
\end{multline}
as $n\to+\infty$. 

Let us now fix a sequence $\eta_n\to 0^+$ such that
\begin{equation}\label{eta}
\max\{\omega_n,\ln|\ln\varepsilon_n|\}=o(\eta_n|\ln\varepsilon_n|)\quad\text{as $n\to+\infty$}\,.
\end{equation}
We group the vortex balls into the following classes and we define:
\begin{align*}
&  D_n^*:= \sum_{i\in I_n^*} d_{i,n}\,, \qquad
I_n^*:=\big\{i\in I_n\,: \; d_{i,n}\geq 0 \;\mbox{and } \,  \zeta(a^n_i)>\zeta_{\rm max}-\eta_n \big\}\,; \\
& D_n^+= \sum_{i\in I_n^+} d_{i,n}\,, \qquad
I_n^+:=\big\{i\in I_n\,: \; d_{i,n}\geq 0  \;\mbox{and } \,  \zeta(a^n_i)\leq\zeta_{\rm max}- \eta_n\big\}\,; \\
&  D_n^-:= \sum_{i\in I_n^-} |d_{i,n}|\,, \qquad
 I_n^-:=\big\{i\in I_n\,:\; \ d_{i,n} < 0\big\}\,.
\end{align*}
Observe that $D_n=D_n^*+D_n^+ +D_n^-$. We claim that 
\begin{align}
\label{prebdD+} D_n^+& \leq C \,\frac{\max\{\omega_n,\ln|\ln\varepsilon_n|\}D_n}{ \eta_n|\ln\varepsilon_n|}  \,, \\[5pt]
\label{prebdD-} D_n^- &\le C \, \frac{\max\{\omega_n,\ln|\ln\varepsilon_n|\}D_n}{|\ln\varepsilon_n|} \,, \\
\label{prebdDn} D_n &\le  C  \max\{\omega_n,\ln|\ln\varepsilon_n|\}\,,
\end{align}
for a constant $C>0$ independent of $n$.  In particular, if $\omega^2_n=o(|\ln\varepsilon_n|)$, then
we can choose $\eta_n$ satisfying in addition $\max\{\omega^2_n,(\ln|\ln\varepsilon_n|)^2\}=o(\eta_n|\ln\varepsilon_n|)$, and consequently  \eqref{prebdDn} yields 
$D_n^+=D_n^- = 0$ for $n$ large enough. 
\vskip5pt

We evaluate the lower bound for each class of vortex ball separately.  First, we use the explicit form of $\Omega_n$ (see (\ref{omega})) and the bound 
$\zeta(x)\le\zeta_{\rm max}$ to obtain,
\begin{equation}\label{lb*}
\pi D_n^*|\ln\varepsilon_n| - 2\pi \Omega_n \sum_{i\in I_n^*} d_{i,n}\zeta(a^n_i)
  \geq  -2\pi \omega_n \zeta_{\rm max} D_n^*\,.
\end{equation}
For negative degrees we have the simple estimate
\begin{equation}\label{lb-}
 \pi D_n^- |\ln\varepsilon_n|
     - 2\pi\Omega_n\sum_{i\in I_n^-} d_{i,n} \zeta(a^n_i)
\geq \pi D_n^- |\ln\varepsilon_n|\,.
\end{equation}
Then for the vortex balls staying away from $\Sigma$, we have
\begin{multline}\label{lb+}
\pi D_n^+|\ln\varepsilon_n| - 2\pi \Omega_n \sum_{i\in I_n^+} d_{i,n}\zeta(a^n_i) 
\geq   \left( \pi|\ln\varepsilon_n| - 2\pi \Omega_n\zeta_{\rm max}\right) D_n^+\,+\\
         + 2\pi\Omega_n \sum_{i\in I_n^+} d_{i,n} (\zeta_{\rm max}-\zeta(a^n_i))  
\geq D_n^+\left( -2\pi\omega_n\zeta_{\rm max} + 2\pi \Omega_n\eta_n\right) \,\geq \\
 \geq  C\Omega_n\eta_n D_n^+\,,
\end{multline}
since $\omega_n=o(\eta_n\Omega_n)$. 
We now insert  \eqref{lb*}, \eqref{lb-} and \eqref{lb+} into (\ref{lb1}),
\begin{align}
\nonumber  O(\omega_n^2) &\geq \pi D_n \left(|\ln\varepsilon_n|-C_1\ln|\ln\varepsilon_n|\right)
     - 2\pi\Omega_n\sum_{i\in I_n} d_{i,n} \zeta(a^n_i) +
      \int_{\DD\setminus\cup B_i^n} |\nabla u_n|^2dx+ o(1) \\
 \label{lb2}
 & \begin{multlined}[10cm] 
 \geq -\pi C_1 D_n\ln|\ln\varepsilon_n| -2\pi \omega_n \zeta_{\rm max} D_n^* + \pi D_n^- |\ln\varepsilon_n|  +C\eta_n\Omega_n D_n^+\,+ \\
 + \int_{\DD\setminus\cup B_i^n} |\nabla u_n|^2 dx\,. 
\end{multlined}
\end{align}
Rearranging all terms we derive
\begin{align*}
D_n^- |\ln\varepsilon_n| + D_n^+\eta_n\Omega_n &\leq C\left( \omega_n D_n^* + \ln|\ln\varepsilon_n| D_n+\omega_n^2\right) \\
 &  \leq C \max\{\omega_n,\ln|\ln\varepsilon_n|\} D_n\,, 
 \end{align*}
which proves \eqref{prebdD+} and \eqref{prebdD-}. 

To prove  \eqref{prebdDn}, we argue as in \cite[pg. 58--60]{AB1} to obtain a constant $C'>0$ (independent of $n$) such that
\begin{equation}\label{outside}   
\int_{\DD\setminus\cup B_i^n} |\nabla u_n|^2dx \geq C' D_n^2\,.
\end{equation}
Accepting \eqref{outside} we  return to the lower bound 
(\ref{lb2}) to deduce 
$$ D_n^2-C \left(\omega_n + \ln|\ln\varepsilon_n|\right) D_n \leq O(\omega_n^2) \,,  $$
so that $D_n\leq C\max\{\omega_n,\ln|\ln\varepsilon_n|\}$ and estimate \eqref{prebdDn}  is established.

It remains to show \eqref{outside}.  To this aim  we identify an annular band lying outside of 
$\Sigma$ and use the fact that the total degree is approximately constant in that band.
Since the boundary $\partial\DD$ is assumed to be smooth, there exists 
$0<\delta_0<\frac{1}{2} \dist(\Sigma,\partial \DD)$ such that the  function
$\varrho(x):=\dist(x,\partial\DD)$ is smooth in $\DD_{\varepsilon_n}^{\delta_0}$ with
$$  \DD^{\delta_0}_{\varepsilon_n}:= \left\{ x\in\DD \,:\, \varepsilon_n<\varrho(x)\leq \delta_0\right\}\,.$$
The level sets 
$$   {\cal C}_t:= \{x\in\DD: \  \varrho(x)=t\}\,, $$
are  smooth and diffeomorphic to $\partial\DD$ for all $t\in [0,\delta_0]$.
Define the set $T_n\subset [0,\delta_0]$ by 
$$  T_n:= \left\{t\in (\varepsilon_n,\delta_0]\,: \; {\cal C}_t\cap \cup_{i\in I_n} B_i^n = \emptyset\right\}\,.$$
By the choice of $r_n$ and claim {\it (i)} in Proposition \ref{consvortexball}, $T_n$ is a finite union of disjoint intervals 
and $\mathcal{L}^1((\varepsilon_n,\delta_0]\setminus T_n)\leq O(|\ln\varepsilon_n|^{-4})$. 
From claim {\it (ii)} in Proposition \ref{consvortexball},
we can define the degree of $u_n$ on ${\cal C}_t$ for every $t\in T_n$, {\it i.e.}, 
$$  D_n(t):= {\rm deg}\bigg(\frac{u_n}{|u_n|}, {\cal C}_t\bigg)=\frac{1}{2\pi}\int_{{\cal C}_t}\frac{u_n}{|u_n|}\wedge\nabla_\tau\bigg(\frac{u_n}{|u_n|}\bigg)\,d{\cal H}^1\,,$$
where $\nabla_\tau$ denotes the tangential derivative along ${\cal C}_t$ oriented counterclockwise.
Setting $I_n(t):=\{i\in I_n\,:\;\varrho(a_i^n)\geq \rho_{i,n}+t\}$ ($\rho_{i,n}$ being the radius of the ball $B_i^n$), we have 
$$D_n(t)=\sum_{i\in I_n(t)}d_{i,n}\,.$$
Using \eqref{prebdD+} and \eqref{prebdD-} we infer that  for $n$ large enough, 
$$ D_n(t) \geq D_n - 2D_n^- - D_n^+ \geq \frac12 D_n
\quad\text{for every $t\in T_n$}\,. $$
Denote $v_n:=u_n/|u_n|$. In view of claim {\it (ii)} in Proposition \ref{consvortexball}, using the Coarea Formula and Jensen Inequality we can estimate  for $n$ large enough, 
\begin{multline*}
\int_{\DD\setminus\cup B_i^n} |\nabla u_n|^2 dx
\geq \int_{\DD_{\varepsilon_n}^{\delta_0}\setminus\cup B_i^n} |u_n|^2|\nabla v_n|^2dx  
\geq \frac{1}{2}\int_{\varepsilon_n}^{\delta_0}\bigg(\int_{{\cal C}_t} |\nabla v_n|^2d{\cal H}^1\bigg)dt  \,\geq \\
\geq\frac{1}{2}\int_{T_n}\bigg(\int_{{\cal C}_t} |v_n\wedge\nabla_\tau v_n|^2d{\cal H}^1\bigg) dt 
\geq  2\pi^2\int_{T_n} \frac{|D_n(t)|^2}{{\cal H}^1({\cal C}_t)}\, dt \geq  C' D_n^2\,,
\end{multline*}
which completes the proof of \eqref{outside}.
 \vskip5pt

\noindent{\it Step 3.} If $\ln|\ln\varepsilon_n|\leq o(\omega_n)$ the conclusion follows from \eqref{prebdD+}, \eqref{prebdD-} and \eqref{prebdDn}. 
If $\omega_n\leq O(\ln|\ln\varepsilon_n|)$ we must refine our lower bound by growing the vortex balls.  First observe that in this regime, 
\eqref{prebdD+} and \eqref{prebdD-} ensures that 
\begin{equation}\label{vandegnegfar}
D_n^-=D_n^+=0
\end{equation} 
for $n$ large so that each ball $B_i^n\subset \DD_{\varepsilon_n}$ carries a nonnegative degree and $D_n=D_n^*$. 

We choose a new radius $s_n:=e^{-\sqrt{\omega_n}}$ and thus $s_n>r_n$ for $n$ large enough.  We now reapply 
Proposition~\ref{consvortexball} with $r=s_n$  
to obtain a new family of larger balls $\{\tilde B_j^n\}_{j\in J_n}$,  each new ball $\tilde B_j^n$ containing 
one or more of the smaller balls $\{B_i^n\}$.  
By claim {\it (ii)} in Proposition~\ref{consvortexball} and \eqref{vandegnegfar} we have 
 $D_{s_n}=D_n$.  Using the lower bound (\ref{balls3}) together with \eqref{lem2} and \eqref{vandegnegfar},  we can argue as in Step 2 to derive 
\begin{align}
\nonumber O(\omega_n^2) \geq \omega_n^2 F_{\varepsilon_n}(u_n) &\geq
 \pi D_n\left(\ln\left(\frac{s_n}{\varepsilon_nD_n}\right)-C_1\right) - 2\pi\Omega_n\zeta_{\rm max}D_n  + C' D_n^2\\
\label{lbrefined} &\geq \pi D_n\left(\ln\left(\frac{s_n}{D_n}\right)-C_1\right) - 2\pi\omega_n\zeta_{\rm max}D_n  + C' D_n^2\,.
\end{align}
Next we distinguish two cases. First assume that $\ln D_n \leq O(\omega_n)$.  In this case, \eqref{lbrefined} yields the inequality
$  D_n^2-C\omega_n D_n\leq O(\omega_n^2) $ 
(with $C>0$ independent of $n$) so that  $D_n\leq O(\omega_n)$ as claimed. 
If  $\omega_n=o(\ln D_n)$,  we  obtain the bound
$ D_n^2-C D_n\ln D_n \leq O(\omega_n^2) $ 
which also yields $D_n\leq O(\omega_n)$, and the proof of Proposition~\ref{degreebound} is complete.
\end{proof} 

We are now ready to prove claim {\it (i)} in Theorem \ref{mainthm}. 

\begin{theorem}\label{compactness} Assume that \emph{(H1)} and \eqref{omega} hold. 
Let $\varepsilon_n\to 0^+$ and  let $\{u_n\}_{n\in\NN}$ be a sequence in $H^1(\DD;\mathbb{C})$  such that $\sup_n F_{\varepsilon_n}(u_n)<+\infty$. 
Then there exist a subsequence (not relabelled) and a nonnegative Radon measure $\mu$ supported by $\Sigma$ such that  
$$\mu_n:=\frac{1}{\omega_n}\,{\rm curl}\,j(u_n) \mathop{\longrightarrow}\limits_{n\to+\infty} \mu\quad\text{strongly in $\big(C_0^{0,1}(\mathcal{D})\big)^*$}\,.$$
\end{theorem}

\begin{proof}
{\it Step 1.} We first assume that $\{u_n\}_{n\in\NN}\subset H^1(\DD;\mathbb{C})\cap C^1(\mathcal{D})$. Using the notations of the previous proof, we apply Proposition~\ref{degreebound} to 
obtain the family of vortex balls  $\{B_i^n\}_{i\in I_n}$. Define the measure 
$$\bar \mu_n:=\frac{2\pi}{\omega_n}\sum_{i\in I_n}d_{i,n}\delta_{a_i^n} \,.$$
Since $D_n\leq O(\omega_n)$ we have $|\bar \mu_n|(\mathcal{D})\leq C$ for a constant $C$ independent of $n$. Therefore, up to a subsequence, 
$\bar \mu_n\rightharpoonup \mu $ as $n\to+\infty$ 
weakly* in the sense of measures on $\mathcal{D}$ for some finite Radon measure $\mu$. We claim that $\mu$ is nonnegative and supported by $\Sigma$. Indeed, decompose $\bar \mu_n$ in 
its Hahn decomposition, {\it i.e.}, write $\bar \mu_n=\bar \mu_n^+-\bar \mu^-_n$ where $\bar \mu_n^+$ and $-\bar \mu_n^-$ are respectively 
the positive and the negative parts of $\bar \mu_n$. Then we have 
$$\bar\mu_n^-(\mathcal{D})=\frac{2\pi D_n^-}{\omega_n}\to 0\quad\text{as $n\to+\infty$}\,,$$
thanks to \eqref{prebdD-}, and the nonnegativity of $\mu$ follows. Now consider the sequence of  sets $\mathcal{V}_n:=\{\zeta(x)\leq \zeta_{\rm max}-\eta_n\}$ where $\eta_n$ is 
given by \eqref{eta}. In view of \eqref{prebdD+}, we have 
$$\bar\mu_n^+(\mathcal{V}_n)=\frac{2\pi D_n^+}{\omega_n} \to 0\,,$$ 
which clearly implies that ${\rm supp}\,\mu\subset\Sigma$. 

By the compact embedding $(C^0_0(\mathcal{D}))^*\hookrightarrow(C_0^{0,1}(\mathcal{D}))^*$, we deduce that $\bar\mu_n\to \mu$ strongly in $(C_0^{0,1}(\mathcal{D}))^*$. On the other hand, 
\eqref{jacobian} yields 
$$\big\|\mu_n -\bar \mu_n\big\|_{(C_0^{01}(\mathcal{D}))^*} \mathop{\longrightarrow}\limits_{n\to+\infty} 0\,,$$
and the conclusion follows. 
\vskip5pt

\noindent{\it Step 2.} We now consider the general case. In view of the strong continuity of the functional $F_{\varepsilon_n}$ under strong $H^1$-convergence,  
we can find a sequence $\{\tilde u_n\}_{n\in\mathbb{N}}\subset  H^1(\DD;\mathbb{C})\cap C^1(\mathcal{D})$ such that for every $n$, 
\begin{equation}\label{app}
\|u_n-\tilde u_n\|_{H^1(\mathcal{D})}\leq \varepsilon_n\,, 
\end{equation}
and
\begin{equation}\label{enapp}
 F_{\varepsilon_n}(\tilde u_n)\leq F_{\varepsilon_n}(u_n)+1\,,
\end{equation}
so that $\sup_n F_{\varepsilon_n}(\tilde u_n)<+\infty$.
 
Given an arbitrary $\varphi \in C^{0,1}_0(\mathcal{D})$ satisfying $|\nabla\varphi|\leq 1$, we estimate
\begin{multline*}
\bigg|\int_{\mathcal{D}}\big(j(u_n)-j(\tilde u_n)\big)\cdot \nabla^\perp\varphi \,dx\bigg| \leq \,\\
\leq \|u_n-\tilde u_n\|_{L^2(\mathcal{D})}\|\nabla u_n\|_{L^2(\mathcal{D})}+
\|\tilde u_n\|_{L^2(\mathcal{D})}\|\nabla(u_n-\tilde u_n)\|_{L^2(\mathcal{D})}\leq C\varepsilon_n|\ln\varepsilon_n|\,,
\end{multline*}
using \eqref{app} and \eqref{enapp} together with Lemma \ref{lem1}. 
As a consequence, setting $\tilde \mu_n:=\omega_n^{-1}\curl j(\tilde u_n)$, we have
\begin{equation}\label{convapp}
 \big\|\mu_n -\tilde  \mu_n\big\|_{(C_0^{01}(\mathcal{D}))^*} \mathop{\longrightarrow}\limits_{n\to+\infty} 0\,.
 \end{equation}
Applying Step 1 to  $\{\tilde u_n\}$, up to a subsequence we have $\tilde  \mu_n\to \mu$ in $(C_0^{0,1}(\mathcal{D}))^*$ for some 
nonnegative Radon measure $\mu$ supported by $\Sigma$. Then \eqref{convapp} yields $\mu_n\to \mu$ in $(C_0^{0,1}(\mathcal{D}))^*$ 
and the proof is complete.
\end{proof}

\section{The lower bound inequality}                                

This section is devoted to the proof of claim {\it (ii)} in Theorem \ref{mainthm} that we summarize in the following result.

\begin{theorem}\label{lwbd} Assume that \emph{(H1)} and \eqref{omega} hold. 
Let $\varepsilon_n\to0^+$ be an arbitrary sequence and let $\{u_n\}_{n\in\NN}\subset H^1(\mathcal{D};\mathbb{C})$ be such that 
\begin{equation}\label{hypconv}
\mu_n=\frac{1}{\omega_n}\curl j(u_n) \mathop{\longrightarrow}\limits_{n\to+\infty} \mu \quad\text{strongly in $(C^{0,1}_0(\mathcal{D}))^*$}\,,
\end{equation}
for some nonnegative Radon measure $\mu$ supported by $\Sigma$. Then, 
\begin{equation}\label{liminf}
\liminf_{n\to+\infty} F_{\varepsilon_n}(u_n) \geq I(\mu)-\zeta_{\rm max}\,\mu(\mathcal{D})\,.
\end{equation}
In particular, if the left hand side in (\ref{liminf}) is finite, then $\mu\in H^{-1}(\mathcal{D})$. 
\end{theorem}

\begin{proof} We will use in this proof the notations of the previous section. Without loss of generality, we may assume that 
\begin{equation}\label{liminflim}
\liminf_{n\to+\infty} F_{\varepsilon_n}(u_n)=\lim_{n\to+\infty} F_{\varepsilon_n}(u_n)<+\infty\,.
\end{equation}
Moreover, by Step 2 in the proof of Theorem \ref{compactness}, we may also assume that 
$\{u_n\}_{n\in\NN}\subset H^1(\mathcal{D};\mathbb{C})\cap C^1(\mathcal{D})$. We shall distinguish two cases.
\vskip5pt

\noindent{\it Case 1.}  We first assume that $\ln|\ln\varepsilon_n|\leq o(\omega_n)$. We consider the family of vortex balls  $\{B_i^n\}_{i\in I_n}$ constructed in the proof of Proposition~\ref{degreebound}, 
and we refer to it for the notations. Arguing as in \eqref{lem2} we obtain
\begin{align}
\nonumber \Omega_n\int_{\mathcal{D}}\nabla^\perp\zeta\cdot j(u_n)\,dx &=-\frac{\pi |\ln\varepsilon_n|}{\zeta_{\rm max}}\sum_{i\in I_n}d_{i,n}\zeta(a_i^n)
+\omega_n\int_{\mathcal{D}}\nabla^\perp\zeta\cdot j(u_n)\,dx+o(1)\\
\label{newrot} &\geq -\pi D_n |\ln\varepsilon_n|
-\omega^2_n\langle \mu_n,\zeta\rangle+o(1)
\end{align}
as $n\to+\infty$, where $\langle\cdot,\cdot\rangle$ denotes the duality pairing $(C_0^{0,1})^*$--$\,C_0^{0,1}$.

Combining the lower bound \eqref{balls3} with \eqref{newrot}, we infer that 
\begin{align*}
\omega^2_n F_{\varepsilon_n}(u_n) &\geq  \pi D_n\left( \ln\left(\frac{r_n}{\varepsilon_n D_n}\right)-C_1\right)
+ \frac12\int_{\DD\setminus\cup B_i^n} |\nabla u_n|^2 dx
         +
          \Omega_n\int_{\mathcal{D}}\nabla^\perp\zeta\cdot j(u_n)\,dx \\
     & \geq   \pi D_n\left( \ln\left(\frac{r_n}{D_n}\right)-C_1\right)
         + \frac12 \int_{\DD\setminus\cup B_i^n}  |\nabla u_n|^2 dx  - \omega^2_n\langle \mu_n,\zeta\rangle+o(1)\,,
\end{align*}
Since $r_n=|\ln\varepsilon_n|^{-4}$ and $D_n\leq O(\omega_n)$ by Proposition~\ref{degreebound}, dividing the previous inequality by $\omega_n^2$ yields
\begin{align*}
F_{\varepsilon_n}(u_n)&\geq \frac{1}{2\omega_n^2} \int_{\DD\setminus\cup B_i^n}  |\nabla u_n|^2 dx  - \langle \mu_n,\zeta\rangle+o(1)\\
&\geq  \frac{1}{2\omega_n^2} \int_{\DD\setminus\cup B_i^n}  |\nabla u_n|^2 dx  - \zeta_{\rm max}\mu(\mathcal{D}) +o(1)\,.
\end{align*}
In the last inequality, we have used \eqref{hypconv} and the fact that $\mu$ is supported by $\Sigma$.  
In view of claim {\it (ii)} in Proposition \ref{consvortexball}, we estimate
$$ \int_{\DD\setminus\cup B_i^n}  |\nabla u_n|^2 dx\geq \frac{1}{1+\varepsilon_n^{1/4}} \int_{\DD_{\varepsilon_n}\setminus\cup B_i^n}  |j(u_n)|^2 dx\,. $$
Next we define
\begin{equation}\label{deftildej}
 \tilde j_n(x):= 
\begin{cases} 
\omega_n^{-1}j(u_n(x)) & \text{if $x\in\DD_{\varepsilon_n}\setminus\bigcup_{i\in I_n} B_i^n$}\,, \\
0 & \text{otherwise}\,,
\end{cases}
\end{equation}
so that by \eqref{liminflim}, 
\begin{equation}\label{saispa}
O(1)\geq F_{\varepsilon_n}(u_n)\geq  \frac{1}{2} \int_{\DD}  |\tilde j_n(x)|^2 dx  - \zeta_{\rm max}\mu(\mathcal{D}) +o(1)\,.
\end{equation}                                   
Hence there exist  a subsequence $\varepsilon_n\to 0$ (not relabelled) and $j_*\in L^2(\DD;\RR^2)$ such that 
$\tilde j_n\rightharpoonup  j_*$  weakly in $L^2(\DD;\RR^2)$ as $n\to+\infty$.  By lower semicontinuity, we have 
\begin{equation}\label{liminfj}
 \lim_{n\to+\infty} F_{\varepsilon_n}(u_n) \geq \frac{1}{2}\int_\DD  |j_*|^2 dx-  \zeta_{\rm max}\mu(\mathcal{D})\,.  
 \end{equation}
It remains to tie the limit $j_*$ to the limit $\mu$ of the normalized weak Jacobians. 
To this aim we fix  $\varphi\in\mathscr{D}(\DD)$.  
Using Lemma \ref{lem1}, claim {\it (i)} in  Proposition \ref{consvortexball} and H\"older Inequality, we estimate
\begin{align*}
\bigg| \int_{\mathcal\cup B_i^n} \nabla^\perp\varphi\cdot j(u_n)\, dx\bigg|
     &\leq \big(\mathcal{L}^2\big(\cup_{i\in I_{n}} B_i^n\big)\big)^{1/4} \|\nabla\varphi\|_\infty\|u_n\|_{L^4(\mathcal{D})} \|\nabla u_n\|_{L^2(\mathcal{D})}\\
        & \le  C\, r_n^{1/2}\, |\ln\varepsilon_n |=o(1)\,. 
\end{align*}
Since ${\rm supp}\,\varphi\,\subset \mathcal{D}_{\varepsilon_n}$ for $n$ large enough, we  deduce that 
\begin{multline*}
\int_\DD  \nabla^\perp\varphi\cdot j_*\, dx =
\lim_{n\to +\infty}  \int_\DD \nabla^\perp\varphi\cdot \tilde j_n\, dx 
= \lim_{n\to +\infty} \omega_n^{-1}\int_{\DD_{\varepsilon_n}\setminus\cup B_i^n}  \nabla^\perp\varphi\cdot j(u_n)\, dx \,= \\
= \lim_{n\to+\infty}  \omega_n^{-1}\int_{\DD} \nabla^\perp\varphi\cdot j(u_n)\, dx  
=-\lim_{n\to+\infty}\,\langle\mu_n,\varphi\rangle  
= - \int_\DD  \varphi\, d\mu\,.
\end{multline*}
Consequently, 
\begin{equation}\label{curlj}
\curl j_*=\mu\quad\text{in $\mathscr{D}'(\mathcal{D})$}\,.
\end{equation}
In particular, $\mu\in H^{-1}(\mathcal{D})$ since $j_*\in L^2(\mathcal{D};\mathbb{R}^2)$. 

Next we introduce $h_\mu\in H_0^1(\mathcal{D})$ to be  the unique solution of 
$$\begin{cases} 
-\Delta h_\mu=\mu &\text{in $H^{-1}(\mathcal{D})$}\,,\\
h_\mu=0 & \text{on $\partial\mathcal{D}$}\,.
\end{cases}
$$
In view of \eqref{curlj} and the definition of $h_\mu$, we have 
$${\rm curl}\big(j_*+\nabla^\perp h_\mu\big)=0 \quad\text{in $H^{-1}(\mathcal{D})$}\,,$$
so that we can find $g_\mu\in H^1(\mathcal{D})$ satisfying $\nabla g_\mu=j_*+\nabla^\perp h_\mu$. Therefore,  
$$\int_{\mathcal{D}}|j_*|^2dx=\int_{\mathcal{D}}|\nabla h_\mu|^2dx+\int_{\mathcal{D}}|\nabla g_\mu|^2dx \,,$$
since an integration by parts yields $\int_\DD \nabla^\perp h_\mu \cdot\nabla g_\mu =0$ (using the fact $h_\mu$ is constant on $\partial \DD$). 
Going back to \eqref{liminfj}, we infer that
$$ \lim_{n\to+\infty} F_{\varepsilon_n}(u_n) \geq \frac{1}{2}\int_\DD  |\nabla h_\mu|^2 dx- \zeta_{\rm max}\mu(\mathcal{D})\,. $$
On the other hand, using the Green representation of $h_\mu$, we have
$$ \frac{1}{2} \int_\DD |\grad h_\mu|^2\, dx =\frac{1}{2} \iint_{\DD\times\DD}  G(x,y)\, d\mu(x)\, d\mu(y)=I(\mu)\,,
$$
and the conclusion follows. 
\vskip5pt 

\noindent{\it Case 2.}  We now treat the case $\omega_n\leq O(\ln|\ln\varepsilon_n|)$. Consider the family of vortex 
balls $\{\tilde B_j^n\}_{j\in J_n}$ (of size $s_n=e^{-\sqrt{\omega_n}}$) constructed in Step 3 in the proof of Proposition \ref{degreebound}. 
Recall that this family satisfies $D_{s_n}=D_n$ for $n$ large. Combining the lower bound \eqref{balls3} for the family  $\{\tilde B_j^n\}_{j\in J_n}$ with \eqref{newrot}, 
we derive  
\begin{align*}
\omega^2_n F_{\varepsilon_n}(u_n) &\geq  \pi D_n\left( \ln\left(\frac{s_n}{D_n}\right)-C_1\right)
         + \frac12 \int_{\DD\setminus\cup\tilde B_j^n}  |\nabla u_n|^2 dx  - \omega^2_n\langle \mu_n,\zeta\rangle+o(1)\,.
\end{align*}
Then arguing as in \eqref{saispa}, we infer that
\begin{equation*}
F_{\varepsilon_n}(u_n)\geq  \frac{1}{2} \int_{\DD}  |\hat j_n(x)|^2 dx  -\zeta_{\rm max}\mu(\mathcal{D}) +o(1)\,,
\end{equation*}  
where
$$\hat j_n(x):=\begin{cases}
\omega_n^{-1}j(u_n(x)) & \text{if $x\in \mathcal{D}_{\varepsilon_n}\setminus \cup_{j\in J_n}\tilde B_j^n$}\,,\\
0 & \text{otherwise}\,.
\end{cases}
$$
As previously, up to a subsequence we have $\hat j_n\rightharpoonup j_*$ weakly in $L^2(\mathcal{D};\mathbb{R}^2)$, and 
$$\lim_{n\to+\infty} F_{\varepsilon_n}(u_n)\geq  \frac{1}{2} \int_{\DD}  |j_*(x)|^2 dx  -  \zeta_{\rm max}\mu(\mathcal{D})\,. $$
Now it remains to show that $\curl j_*=\mu$ in $\mathscr{D}'(\mathcal{D})$, and then  the proof can be completed as in Step 1. 

We proceed as before, taking an arbitrary $\varphi\in \mathscr{D}(\mathcal{D})$
and using the weak formulation of the Jacobians.   
The key observation is that the contribution of the vortex balls will be negligible provided we can restrict our choice  
of  test functions $\varphi$ to functions {\em constant} in each vortex ball.  This can be achieved thanks to \cite[Proposition~9.6]{SSbouquin}, {\it i.e.},  
given an arbitrary $\varphi\in \mathscr{D}(\mathcal{D})$, there exists a modified function $\tilde\varphi_n$ which is constant on each ball $B_j^n$ and such that
\begin{equation}\label{approx}
 \|\varphi-\tilde\varphi_n\|_{C^{0,\alpha}(\mathcal{D})} \leq Cs_n^{1-\alpha}\,, \qquad
    \|\nabla\varphi-\nabla\tilde\varphi_n\|_{L^1(\mathcal{D})}\le Cs_n
\end{equation}
for each $0\leq\alpha\leq 1$.  Moreover, $\tilde\varphi_n$ has compact support in $\mathcal{D}$ for $n$ large enough. 
From \eqref{approx} and the $L^2$-boundedness of the normalized currents $\hat j_n$, we derive
\begin{equation}\label{approx2}
 \left|\int_\DD\left(\nabla^\perp\varphi - \nabla^\perp\tilde\varphi_n\right)\cdot  \hat j_n \,dx\right| \leq  
 \|\hat j_n\|_{L^2(\mathcal{D})} \|\nabla\varphi-\nabla\tilde\varphi_n\|_{L^2(\mathcal{D})}
\leq C s_n^{1/2}\,.  
\end{equation}
Then \eqref{approx}, \eqref{approx2} and the strong convergence of $\mu_n$ to $\mu$ yield
\begin{multline*}
\int_\DD \nabla^\perp\varphi\cdot  j_*\, dx =
\lim_{n\to +\infty}  \int_\DD \nabla^\perp\varphi\cdot\hat j_n\, dx  
= \lim_{n\to +\infty}  \int_\DD \nabla^\perp\tilde\varphi_n\cdot\hat j_n\, dx \,= \\
=  \lim_{n\to +\infty} \omega_n^{-1} \int_{\DD}   \nabla^\perp\tilde\varphi_n\cdot j(u_n)\, dx  =
\lim_{n\to+\infty}\bigg( -\int_{\mathcal{D}}\tilde\varphi_n d\mu +\langle \mu-\mu_n, \tilde\varphi_n\rangle\bigg)
= - \int_\DD  \varphi\, d\mu\,,
\end{multline*}
and the conclusion follows.
\end{proof}

\section{The upper bound inequality}\label{secupbd}                  

Throughout this section we shall use the following notation. For a nonnegative Radon measure $\mu\in H^{-1}(\DD)$ compactly supported in $\DD$, 
we denote by $h_{\mu}\in H^1_0(\DD)$ the solution of 
\begin{equation}\label{defhmu}
\begin{cases} 
-\Delta h_\mu=\mu &\text{in $H^{-1}(\mathcal{D})$}\,,\\
h_\mu=0 & \text{on $\partial\mathcal{D}$}\,.
\end{cases}
\end{equation}
Using the Green representation $h_\mu(x) =\int_{\DD}G(x,y)\,d\mu(y)$, we have 
$$I(\mu)=\frac{1}{2}\int_{\DD}|\nabla h_\mu|^2\,dx\,,$$
and the following elementary lemma holds.

\begin{lemma}\label{convIgen}
Let $\{\mu_n\}_{n\in\NN}$ be a sequence of  nonnegative Radon measure in $H^{-1}(\DD)$ with compact support in $\DD$. Assume that $\mu_n\rightharpoonup\mu$ 
weakly* as measures on $\DD$ as $n\to+\infty$, for some $\mu\in H^{-1}(\DD)$ compactly supported in $\DD$. Then $I(\mu_n)\to I(\mu)$ if and only if $h_{\mu_n}\to h_\mu$ strongly in $H^1(\DD)$ as $n\to+\infty$.  
\end{lemma}

We may now start the proof of claim {\it (iii)} in Theorem \ref{mainthm} in the case where the measure $\mu$ is absolutely continuous with respect to $\mathcal{H}^1\res\Sigma$. 

\begin{proposition}\label{limsupabscont}
Assume that \emph{(H1)}, \emph{(H2)} and \eqref{omega} hold. 
Let $\varepsilon_n\to0^+$ be an arbitrary sequence. For every nonnegative Radon measure $\mu$ of the form 
\begin{equation}\label{abscontmeas}
\mu=f(x)\,\mathcal{H}^1\res \Sigma
\end{equation}
with $f\in L^\infty(\Sigma)$, there exists a sequence $\{u_n\}_{n\in\mathbb{N}}\subset H^1(\mathcal{D};\mathbb{C})$ such that
$$\frac{1}{\omega_n}\,\curl\,j(u_n)\mathop{\longrightarrow}\limits_{n\to+\infty} \mu\quad\text{strongly in $(C_0^{0,1}(\mathcal{D}))^*$}\,, $$
and
$$\lim_{n\to+\infty}\,F_{\varepsilon_n}(u_n)= I(\mu)- \zeta_{\rm max}\mu(\mathcal{D})\,. $$
\end{proposition}

\begin{proof}  
Without loss of generality we may assume that $f\not\equiv 0$, the case $\mu=0$ being easily true. 
It is well known that a measure of the form \eqref{abscontmeas}
 belongs to $H^{-1}(\mathcal{D})$, see {\it e.g.} \cite[Theorem~4.7.5]{Zi}.

We recall that the Dirichlet Green's function $G$ in $\DD$ defined by \eqref{defGreen} satisfies 
\begin{enumerate}
\item[\it (i)] $G(x,y)\ge 0$  for every 
$x\in \overline{\DD}\setminus\{y\}$  and for every $y\in\mathcal{D}$; 

\item[\it (ii)] for any compact set $K\subset\subset\DD$ there exists a constant $C_K$
such that
\begin{equation}\label{Green} 
\left| G(x,y)+\frac{1}{2\pi} \ln |x-y|
\right|\leq C_K
\end{equation}
for all $y\in K$ and $x\in \overline{\DD}$.
\end{enumerate}

In the three first steps below, we assume that the density function $f$ does not vanish on $\Sigma$, {\it i.e.}, $f\geq \delta$ $\mathcal{H}^1$-a.e. on $\Sigma$ 
for some constant $\delta>0$. The general case is considered in Step~4. 
\vskip5pt

\noindent{\it Step 1.} 
We will construct a trial function using the Green function $G(x,y)$  in the spirit of \cite{SS3}. 
Let $\gamma: [0,\ell]\to\Sigma$ be an arclength parametrization of the curve
$\Sigma$ (so that $\ell=\mathcal{H}^1(\Sigma)$), and define for $t\in[0,\ell]$, 
$$  M(t):=\mu\big(\gamma([0,t])\big)\,. $$
Then $M(\cdot)$ is strictly increasing, $M(0)=0$, and  $M(\ell)=\mu(\Sigma)=\mu(\mathcal{D})$.  Moreover $M(\cdot)$ is continuous since $\mu\in H^{-1}(\DD)$ and thus atomless. 

Next we introduce for $n$ large enough,
$$ D_n:=\left[\frac{\omega_n\mu(\Sigma)}{2\pi}\right]\,,$$
where $[\cdot]$ denotes the integer part. 
Since $M$ is continuous and increasing, we can define for $k=0,\ldots, D_n$,
$$t_{k,n}:=M^{-1}\big(2\pi k\omega_n^{-1}\big)\,.$$
Now we set for $k=0,\ldots,D_n$,
$$a^n_{k}:=\gamma(t_{k,n})\,,$$
and we claim that 
\begin{equation}\label{space}
C \omega_n^{-1}\leq |a^n_{k}-a^n_{k-1}|\leq 2\pi\delta^{-1} \omega_n^{-1}
\end{equation}
for each $k\in\{1,\ldots,D_n\}$ and for some constant $C>0$ independent of $n$. Write $\Sigma_k:=\gamma([t_{k-1,n},t_{k,n}])$ for $k=1,\ldots,D_n$ so that 
$\Sigma_k$ is a smooth curve whose end-points are $a^n_{k}$ and $a^n_{k-1}$. 
By construction, we have 
$$
2\pi \omega_n^{-1}=M(t_{k,n})-M(t_{k-1,n})=\mu\big(\Sigma_k\big)
\geq\delta \mathcal{H}^1\big(\Sigma_k\big) \geq \delta |a^n_{k}-a^n_{k-1}|\,.
$$
Now the curve $\Sigma$ being smooth and  $|a^n_{k}-a^n_{k-1}|$ small for $n$ large, we deduce
$$ 2\pi\omega_n^{-1}=\mu\big(\Sigma_k\big)\leq \|f\|_{L^\infty(\Sigma)} \mathcal{H}^1\big(\Sigma_k\big) \leq  
C \|f\|_{L^\infty(\Sigma)} |a^n_{k}-a^n_{k-1}| \,.$$

Define the family of measures 
$$  \bar\mu_n := \frac{2\pi}{\omega_n}\sum_{k=1}^{D_n} \delta_{a^n_k}\,. $$
We claim that $\bar\mu_n\rightharpoonup \mu$ weakly* as measures on $\DD$ and $\bar\mu_n(\mathcal{D})\to \mu(\mathcal{D})$. 
To prove the weak* convergence of $\bar\mu_n$, we fix an  
arbitrary  function $\psi\in C^0_0(\mathcal{D})$. Observe that by \eqref{space} and the smoothness of 
$\Sigma$, we have ${\rm diam}(\Sigma_k)\leq C_0\omega_n^{-1}$ for a constant $C_0$  
independent of $k$ and $n$. Therefore,  using $\mu(\Sigma_k)=2\pi \omega_n^{-1}$ we derive
$$
\bigg|\int_\DD  \psi\, d\mu - \int_\DD \psi\, d\bar\mu_n\bigg|
= \bigg| \sum_{k=1}^{D_n}  \int_{\Sigma_k}  \big(\psi(x)-\psi(a_k^n)\big)d\mu\bigg|+o(1)\leq C \mu(\mathcal{D}) {\rm osc}(\psi,C_0\omega_n^{-1})+o(1)\,,
$$
where 
$${\rm osc}(\psi,C_0\omega_n^{-1}):=\sup_{|x-y|\leq C_0\omega_n^{-1}} |\psi(x)-\psi(y)|\;\mathop{\longrightarrow}\limits_{n\to+\infty}0\,,$$ 
and the claim is proved.

Next we must regularize the measure $\bar\mu_n$.  Let us define for $ k=1,\dots,D_n$, 
$$  f_{k,n}(x):=\frac{1}{\pi \varepsilon_n^2}\,\chi_{B_{\varepsilon_n}(a_k^n)}\,,$$
where $\chi_{{B_{\varepsilon_n}(a_k^n)}}$ denotes the characteristic
function of the ball $B_{\varepsilon_n}(a_k^n)$. Set
\begin{equation}\label{deffn}
 f_n(x): =\frac{2\pi}{\omega_n} \sum_{k=1}^{D_n} f_{k,n}(x)\quad\text{and}  \quad \hat \mu_n:=f_n(x)\mathcal{L}^2\res\mathcal{D}\,.
 \end{equation}
Since $\varepsilon_n=o(\omega_n^{-1})$, the functions  $f_{k,n}$ have disjoint supports for $n$ large by \eqref{space}. As a consequence, 
$\hat \mu_n(\mathcal{D})=\bar\mu_n(\mathcal{D})\to \mu(\mathcal{D})$. Since $\bar\mu_n\rightharpoonup\mu$, 
one may also easily check that $\hat\mu_n\rightharpoonup\mu$  weakly* as measures on $\DD$.
\vskip5pt

\noindent{\it Step 2.} According to \eqref{defhmu}, we introduce
$$ h_n:=h_{\omega_n\hat\mu_n} =\omega_n h_{\hat\mu_n}\,.$$ 
Then
\begin{align}
\label{U1}  \int_\DD  |\nabla h_n|^2 dx 
        &= \omega^2_n\iint_{\DD\times\DD} G(x,y)\, f_n(x)f_n(y)\, dxdy \\
        \label{U2.5} & =  4\pi^2\sum_{i,j=1}^{D_n}  \iint_{\DD\times\DD}
       G(x,y) f_{i,n}(x) f_{j,n}(y)\, dx dy \,.
\end{align}
We need to estimate the integral term in the right handside of \eqref{U1}.  
We proceed as in \cite{AAB} and we provide some details
for the reader's convenience. Let ${\cal N}_0\subset\!\subset\DD$ be a small tubular neighborhood of $\Sigma$. 
 Let $0<\alpha<1$ be given small, and set $\Delta_\alpha:=\{ (x,y)\in\DD\times\DD\,: \;
|x-y|<\alpha\}$.  Since $G$ is continuous on $(\mathcal{N}_0\times
\mathcal{N}_0)\setminus \Delta_\alpha$, and the support of  $\hat \mu_n$ lies in $\mathcal{N}_0$ for $n$ large,
by the weak* convergence of $\hat\mu_n$ to $\mu$  we have
\begin{multline}\label{U3}
I_\alpha:=\lim_{n\to+\infty}\iint_{(\DD\times
\DD)\setminus\Delta_\alpha}
G(x,y) \, f_n(x) f_n(y)dx dy\,=\\
=\iint_{(\mathcal{N}_0\times \mathcal{N}_0)\setminus\Delta_\alpha} G(x,y) d\mu(x) d\mu(y) 
\leq  \iint_{\DD\times \DD} G(x,y) d\mu(x) d\mu(y)= 2 I(\mu)\, .
\end{multline}
 Near the diagonal $\Delta_\alpha$ we split the sum in
(\ref{U2.5}) in two terms. Using (\ref{Green}) we estimate
\begin{align}
\nonumber II^n_\alpha:=4\pi^2\sum_{i=1}^{D_n} \iint_{\Delta_\alpha}
G(x,y)\,& f_{i,n}(x) f_{i,n}(y)dxdy =\\
\nonumber &=4 \sum_{i=1}^{D_n} \iint_{B_{1}(0)\times B_{1}(0)} 
    G\left(  a_i^n + \varepsilon_nz_1, a_i^n+\varepsilon_n z_2\right) dz_1 dz_2  \\
\nonumber&\leq  4
    \sum_{i=1}^{D_n} \iint_{B_1(0)\times B_1(0)} 
\left(  \frac{1}{2\pi}
         \ln\left( \frac{1}{\varepsilon_n |z_1-z_2|}\right) + C
     \right) dz_1 dz_2 \\
\label{sumdiag} &\leq 2\pi D_n |\ln\eps_n| + O(\omega_n)\,,
\end{align}
and
\bigskip
\begin{multline}\label{diagdifpt}
III^n_\alpha:=4\pi^2\!\!\!\sum_{0<|a^n_i-a^n_j|<\alpha} 
\iint_{\Delta_\alpha}
 G(x,y)\, f_{i,n}(x)\, f_{j,n}(y)dxdy \,\leq\\
 \leq   C\!\!\!
\sum_{0<|a^n_i-a^n_j|<\alpha}  \left|
\ln|a^n_i-a^n_j|\right| \,.
\end{multline}
By the smoothness of $\Sigma$ and \eqref{space}, there exists a constant $c_0>0$ independent of $n$ such that for every $i\not=j$ and every $(x,y)\in \Sigma_i\times\Sigma_j$, 
$$|a_i^n-a_j^n|\geq c_0|x-y|\,. $$
Since $\mu\in H^{-1}(\mathcal{D})$ and it is supported by $\Sigma$, the map $(x,y)\mapsto \big|\ln (c_0|x-y|)\big|$ belongs to $L^1(\Sigma\times\Sigma,\mu\otimes\mu)$. Therefore, 
by the Mean Value Theorem, for every $i\not=j$ we can find  a pair $(x_i^n,y_j^n)\in \Sigma_i\times\Sigma_j$ such that 
$$\frac{\omega_n^2}{4\pi^2}\iint_{\Sigma_i\times\Sigma_j}\big|\ln (c_0|x-y|)\big|d\mu(y)d\mu(x)= \big|\ln(c_0|x_i^n-y_j^n|)\big|\,,$$
noticing that $4\pi^2\omega_n^{-2}=\mu(\Sigma_i)\mu(\Sigma_j)$. 
Applying the previous inequality to $(x_i^n,y_j^n)$ we deduce from \eqref{diagdifpt} that for $n$ large, 
\begin{align}
\nonumber III^n_\alpha&\leq C \sum_{0<|a^n_i-a^n_j|<\alpha}  \big|\ln(c_0|x_i^n-y_j^n|)\big|\\
\nonumber&\leq C\omega_n^2 \sum_{0<|a^n_i-a^n_j|<\alpha}\iint_{\Sigma_i\times\Sigma_j}\big|\ln (c_0|x-y|)\big|d\mu(y)d\mu(x) \\
\label{diagdifpt2}&\leq C\omega_n^2 \iint_{(\Sigma\times\Sigma)\cap \Delta_{2\alpha}}\big|\ln (c_0|x-y|)\big|d\mu(y)d\mu(x)\,.
\end{align}
On the other hand $\mu\in H^{-1}(\mathcal{D})$ so it is atomless, 
and thus $\mu\otimes \mu$ does not charge $\{x=y\}\cap \DD\times\DD$.  Consequently the integral term in 
the right handside of \eqref{diagdifpt2} vanishes as $\alpha\to 0^+$. Hence,
\begin{equation}\label{diagdifpt3}
\lim_{\alpha\to 0^+}\limsup_{n\to+\infty}\,\omega_n^{-2} III^n_\alpha =0\,.
\end{equation}
Gathering \eqref{U3}, \eqref{sumdiag} and \eqref{diagdifpt3} yields 
\begin{multline*}
\limsup_{n\to+\infty} \bigg(\iint_{\DD\times\DD} G(x,y)\, f_n(x)f_n(y)\, dxdy- 2\pi D_n |\ln\varepsilon_n| \,\omega_n^{-2}\bigg)\leq \\
\leq\limsup_{\alpha\to0^+}\left(I_\alpha+\limsup_{n\to+\infty} \big((II^n_\alpha- 2\pi D_n |\ln\varepsilon_n|)\omega_n^{-2} +\omega_n^{-2}III^n_\alpha\big)\right)\leq 2I(\mu)\,,
\end{multline*}
and in view of \eqref{U1}, we conclude that 
\begin{equation}\label{U4} 
\frac{1}{2}\int_\DD  |\nabla h_n|^2 dx 
\leq    \pi D_n \,|\ln\varepsilon_n|
     + \omega_n^2 I(\mu) + o(\omega_n^2)\,.
\end{equation}
\vskip5pt

\noindent{\it Step 3.}  We shall now define a complex-valued order parameter $u_n$ associated to $h_n$. We proceed as follows. 
Since
\begin{equation}\label{linkcurl}
  \curl\left(-\nabla^\perp h_n\right)
  = -\Delta h_n  = \omega_n f_n  
  \end{equation}
is supported by $\cup_k \overline B_{\varepsilon_n}(a_k^n)$, we may  locally define a phase $\phi_n$ in 
$\DD\setminus\cup_k \overline B_{\varepsilon_n}(a_k^n)$
by
$$   \nabla \phi_n(x) = -\nabla^\perp h_n(x)
    \quad \text{for }x\in \DD\setminus\cup_k \overline B_{\varepsilon_n}(a_k^n)\,.$$
In fact, since the balls $B_{\varepsilon_n}(a_k^n)$ are pairwise disjoint (assuming $n$ large enough) and  the mass of $\omega_n f_n$ is quantized in each such ball, it is easy to show that  
$\phi_n$ is single-valued modulo $2\pi$, {\it i.e.}, for any smooth Jordan curve $\Theta$ 
inside $\DD\setminus\cup_k \overline B_{\varepsilon_n}(a_k^n)$,
$$
\frac{1}{2\pi}\int_\Theta \nabla\phi_n\cdot\tau \, \in 
  \mathbb{Z}\,,
$$
where $\tau:\Theta\to\mathbb{S}^1$ is any  smooth  vector field tangent to $\Theta$.  Hence 
$\exp (i\phi_n(x))$ is well defined for every $x\in\DD\setminus\cup_k \overline B_{\varepsilon_n}(a_k^n)$. 

Then consider a smooth cut-off function $\rho:\mathbb{R}\to [0,1]$ such that  
$\rho(t)\equiv 1$ for $t\geq 2$, and $\rho(t)\equiv 0$ for $t \leq 1$. 
Define 
\begin{equation}\label{defprof}
\rho_n(x) :=
\begin{cases}
\displaystyle \rho\left(\frac{|x-a_k^n|}{\varepsilon_n}\right) & \text{if $x\in B_{2\eps_n}(a_k^n)$ for some  $k=1,\dots,D_n$}\,,\\[10pt] 
1 & \text{otherwise}\,,
\end{cases}
\end{equation}
and observe that 
\begin{equation}\label{enmod}
E_{\varepsilon_n}(\rho_n,\DD)=O(\omega_n) \,.
\end{equation}
Then set
$$u_n(x):=\begin{cases}
\displaystyle \rho_n(x) e^{i\phi_n(x)} & \text{for $x\in \DD\setminus\cup_k \overline B_{\varepsilon_n}(a_k^n)$}\,,\\
0 & \text{otherwise}\,.
\end{cases}
$$ 
One may easily check that $u_n\in H^1(\DD;\mathbb{C})$. We claim that 
\begin{equation}\label{strconvcons}
\mu_n=\omega_n^{-1}\curl\,j(u_n)\mathop{\longrightarrow}\limits_{n\to+\infty} \mu \quad\text{strongly in $(C_0^{0,1}(\DD))^*$}\,.
\end{equation}
A simple computation gives
$$j(u_n)=\rho_n^2\nabla\phi_n =-\rho_n^2 \nabla^\perp h_n\quad \text{a.e. in $\DD$}\,.$$
Given $\varphi\in C_0^{0,1}(\DD)$ satisfying $|\nabla\varphi|\leq 1$, we deduce from \eqref{linkcurl} and the previous indentity, 
\begin{align}
\nonumber \langle \mu_n,\varphi\rangle &=  \frac{1}{\omega_n}\int_\DD\nabla\varphi\cdot\nabla h_ndx+ \frac{1}{\omega_n}\int_{\DD}(\rho_n^2-1)\nabla\varphi\cdot\nabla h_ndx\\
\label{estijac1} & = \langle \hat \mu_n,\varphi\rangle+ \frac{1}{\omega_n}\int_{\cup_k B_{2\eps_n}(a_k^n)}(\rho_n^2-1)\nabla\varphi\cdot\nabla h_ndx\,.
\end{align}
In view of the compact embedding $(C^0_0(\mathcal{D}))^*\hookrightarrow(C_0^{0,1}(\mathcal{D}))^*$,  
$\hat \mu_n\to \mu$ strongly in $(C_0^{0,1}(\mathcal{D}))^*$. Hence we can estimate using \eqref{U4},  
\begin{multline}\label{estijac2}
|\langle \mu_n-\mu,\varphi\rangle | \leq  \big\|\hat \mu_n -  \mu\big\|_{(C_0^{0,1}(\mathcal{D}))^*} +\omega_n^{-1}\big(\mathcal{L}^2(\cup_k B_{2\eps_n}(a_k^n))\big)^{1/2}\|\nabla h_n\|_{L^2(\DD)}\,\leq\\
\leq \big\|\hat \mu_n - \mu\big\|_{(C_0^{0,1}(\mathcal{D}))^*} + C\varepsilon_n|\ln\varepsilon_n|^{1/2}\mathop{\longrightarrow}\limits_{n\to+\infty} 0\,,
\end{multline}
and \eqref{strconvcons} is proved.

We now compute the energy $F_{\varepsilon_n}(u_n)$.  We infer from \eqref{U4} and  \eqref{enmod} that 
 \begin{multline}\label{GLentest}
E_{\varepsilon_n}(u_n,\DD)
  = E_{\varepsilon_n}(\rho_n,\DD)+\frac{1}{2}  \int_{\DD\setminus \cup_k \overline B_{\eps_n}(a_k^n)}  \rho_n^2 |\nabla \phi_n|^2  dx \,\leq \\
  \leq   \frac{1}{2}\int_\DD |\nabla h_n|^2 + O(\omega_n)\leq \pi D_n \,|\ln\varepsilon_n|
     + \omega_n^2 I(\mu) + o(\omega_n^2)\,,
\end{multline}
and it remains to evaluate the interaction with the rotation potential.  First \eqref{strconvcons} yields 
\begin{multline*}  
\Omega_n \int_\DD  \nabla^\perp \zeta \cdot j(u_n)\,dx =\frac{|\ln\varepsilon_n|}{2\zeta_{\rm max}}  \int_\DD  \nabla^\perp \zeta \cdot j(u_n)\,dx - \omega_n^2\langle \mu_n,\zeta\rangle\,=\\
= \frac{|\ln\varepsilon_n|}{2\zeta_{\rm max}}  \int_\DD  \nabla^\perp \zeta \cdot j(u_n)\,dx - \zeta_{\rm max} \mu(\mathcal{D})\omega_n^2+o(\omega_n^2)\,.
  \end{multline*}
Arguing as in \eqref{estijac1}-\eqref{estijac2} we derive 
\begin{multline*}
\int_\DD  \nabla^\perp \zeta \cdot j(u_n)\,dx=-\omega_n\langle\hat\mu_n,\zeta\rangle+O(\varepsilon_n\omega_n|\ln\varepsilon_n|^{1/2})\,= \\
= -2\pi \sum_{k=1}^{D_n} \frac{1}{\pi\varepsilon_n^2}\int_{B_{\varepsilon_n}(a_k^n)}\zeta(x)dx +o(\varepsilon_n|\ln\varepsilon_n|^{3/2}) 
= -2\pi D_n \zeta_{\rm max} +o(\varepsilon_n|\ln\varepsilon_n|^{3/2})\,,
\end{multline*}
and consequently, 
\begin{equation}\label{U5}
\Omega_n \int_\DD  \nabla^\perp \zeta \cdot j(u_n)\,dx = -\pi D_n |\ln\varepsilon_n| - \zeta_{\rm max} \mu(\mathcal{D})\omega_n^2+o(\omega_n^2)\,.
\end{equation}
Combining \eqref{GLentest} with \eqref{U5} finally leads to 
$$F_{\varepsilon_n}(u_n)\leq I(\mu)-\zeta_{\rm max} \mu(\mathcal{D}) +o(1)\,.$$
In view of Theorem \ref{lwbd}, the conclusion follows taking the $\limsup$ as $n\to+\infty$ in the previous inequality. 
\vskip5pt

\noindent{\it Step 4.} We now consider the case where the density $f$ is allowed to vanish. 
Let $\{\delta_k\}\subset \mathbb{R}$ be a sequence decreasing to 0 as $k\to+\infty$. Then for $k\in\mathbb{N}$,  we consider 
the measure 
$$\mu_k:=\mu+\delta_k\mathcal{H}^1\LL \Sigma=\big(f(x)+\delta_k\big)\mathcal{H}^1\LL \Sigma\,.$$ 
By monotone convergence, one has 
\begin{equation}\label{convapproxmeas}
I(\mu_k)\mathop{\longrightarrow}\limits_{k\to+\infty} I(\mu)\qquad\text{and}\qquad \mu_k(\mathcal{D})\mathop{\longrightarrow}\limits_{k\to+\infty} \mu(\mathcal{D})\,.
\end{equation}
Obviously $\mu_k$ also converges to $\mu$ strongly in $(C^{0,1}_0(\mathcal{D}))^*$. 
Applying Step 1 to Step 3,  we find for every $k\in\mathbb{N}$ a sequence $\{v_n^k\}_{n\in\mathbb{N}}\subset H^1(\mathcal{D};\mathbb{C})$ such that 
$\omega_n^{-1}{\rm curl}\,j(v_n^k)\to \mu_k$ strongly in $(C^{0,1}_0(\mathcal{D}))^*$ and $F_{\varepsilon_n}(v_n^k)\to I(\mu_k)-\zeta_{\rm max}\mu_k(\mathcal{D})$ as $n\to+\infty$. 
Hence for every $k\in\mathbb{N}$, we can find $N_k\in\mathbb{N}$ such that for every $n\geq N_k$, 
$$\|\omega_n^{-1}{\rm curl}\,j(v_n^k)-\mu\|_{(C^{0,1}_0(\mathcal{D}))^*}\leq 2^{-k}$$
and
$$\big|F_{\varepsilon_n}(v_n^k)-I(\mu)+\zeta_{\rm max}\mu(\mathcal{D})\big|\leq 2^{-k}\,. $$
Moreover we can assume without loss of generality that the sequence of integers $\{N_k\}_{k\in\mathbb{N}}$ is strictly increasing. Therefore 
given any  integer $n$ large enough, there is a unique $k_n\in\mathbb{N}$ such that $N_{k_n}\leq n < N_{k_n+1}$, and
$k_n\to+\infty$ as $n\to+\infty$.
We may then define $u_n:=  v_{n}^{k_n}$.  By construction,  the sequence $\{u_n\}$ satisfies the required properties. 
\end{proof}

To consider the case of a general measure $\mu$ in $H^{-1}(\DD)$, we shall need the following continuity lemma. 

\begin{lemma}\label{conttrans}
Let $\mu$ be a nonnegative Radon measure in $H^{-1}(\DD)$ such that ${\rm supp}\,\mu\subset\!\subset \DD$. For $\xi\in \RR^2$, let $\tau_\xi\mu$ be the translated measure defined by 
$$\tau_\xi\mu(B)= \mu(-\xi+B)\quad \text{for any Borel set $B\subset \RR^2$}\,.$$
Then there exists $0<\delta<\dist({\rm supp}\,\mu, \partial\DD)$ such that $\tau_\xi\mu\in H^{-1}(\DD)$ for every $\xi\in B_\delta(0)$, and the mapping $\xi\mapsto h_{\tau_\xi\mu}\in H^1(\DD)$
is strongly continuous on $B_\delta(0)$. 
\end{lemma}

\begin{proof}
For $\delta>0$ we set $\widetilde\DD_{\delta}:=\{x\in\RR^2\,,\,\dist(x,\DD)<2\delta\}$. 
Then choose $\delta> 0$ such that $\partial \widetilde\DD_{\delta}$ is smooth and $2\delta<\dist({\rm supp}\,\mu,\partial \DD)$. 
For every $\xi\in B_\delta(0)$ we have $\DD\subset\!\subset \xi+\widetilde\DD_\delta$, ${\rm supp}\,\tau_{\xi}\mu=\xi+{\rm supp}\,\mu\subset \DD$ and 
$\dist(\xi+{\rm supp}\,\mu,\partial \DD)>\delta$. 

Obviously $\mu\in H^{-1}(\widetilde \DD_\delta)$ and we can set 
$\bar h\in H^1(\widetilde\DD_{\delta})$ to be the unique solution of 
$$\begin{cases}
-\Delta \bar h=\mu & \text{in $\widetilde \DD_{\delta}$}\,,\\
 \bar h=0 & \text{on $\widetilde\DD_{\delta}$}\,. 
\end{cases}
$$
By our choice of $\delta$, the function $\bar h$ is smooth in the $\delta$-neighborhood of $\partial\DD$. 
Next, for $\xi\in B_\delta(0)$ we denote by $\bar h_{\xi}\in H^1(\DD)$ the function defined by $\bar h_{\xi}:=\bar h(x-\xi)$ for $x\in \DD$.  
 Observe that $\bar h_\xi\in H^1(\DD)$ and $-\Delta \bar h_\xi= \tau_\xi\mu$ in $\DD$. Hence $\tau_\xi\mu\in H^{-1}(\DD)$. 
 
 Now consider a sequence $\{\xi_n\}\subset B_\delta(0)$ such that $\xi_n\to \xi\in B_\delta(0)$ as $n\to+\infty$. Denote $h_n:=h_{\tau_{\xi_n}\mu}$. We have 
 $\Delta(\bar h_{\xi_n}-h_n)=0$ in $\DD$ and  $(\bar h_{\xi_n}-h_n)=\bar  h(x-\xi_n)$ on $\partial\DD$.  By standard elliptic estimates, $(\bar h_{\xi_n}-h_n)$ strongly converges 
 in $H^1(\DD)$ to the harmonic function in $\DD$ equal to $\bar  h(x-\xi)$ on $\partial\DD$, that is $\bar h_{\xi}-h_{\tau_{\xi}\mu}$. On the other hand, $\bar h_{\xi_n}\to \bar h_{\xi}$ 
 strongly  in $H^1(\DD)$ by strong continuity of translations in $H^1$. Therefore $h_n\to h_{\tau_{\xi}}$ strongly in $H^1(\DD)$, and the proof is complete. 
\end{proof}

\begin{theorem}\label{upbd}
Assume that \emph{(H1)}, \emph{(H2)} and \eqref{omega} hold. Let $\varepsilon_n\to0^+$ be an arbitrary sequence. 
For every nonnegative Radon measure $\mu\in H^{-1}(\mathcal{D})$ 
supported by $\Sigma$, there exists a sequence $\{u_n\}_{n\in\mathbb{N}}\subset H^1(\mathcal{D};\mathbb{R}^2)$ such that
$$\frac{1}{\omega_n}\,\curl\,j(u_n)\mathop{\longrightarrow}\limits_{n\to+\infty} \mu\quad\text{strongly in $(C_0^{0,1}(\mathcal{D}))^*$}\,, $$
and
$$\lim_{n\to+\infty}\,F_{\varepsilon_n}(u_n)= I(\mu)-\zeta_{\rm max}\mu(\mathcal{D})\,. $$
\end{theorem}

\begin{proof} We shall prove that for any nonnegative Radon measure $\mu\in H^{-1}(\mathcal{D})$ 
supported by $\Sigma$, there exists a sequence of nonnegative Radon measures $\{\mu_k\}_{k\in\NN}$ of the form \eqref{abscontmeas} such that $\mu_k\to\mu$ 
strongly in $(C_0^{0,1}(\mathcal{D}))^*$, $\mu_k(\DD)\to\mu(\DD)$ and 
\begin{equation}\label{convapproxpot}
I(\mu_k)\to I(\mu)\quad\text{as $k\to+\infty$}\,.
\end{equation}
Assuming that such a sequence exists, Proposition \ref{limsupabscont} yields for each $k$ a sequence $\{v_n^k\}_{n\in\mathbb{N}}\subset H^1(\mathcal{D};\mathbb{R}^2)$ such that 
$\omega_n^{-1}{\rm curl}\,j(v_n^k)\to \mu_k$ strongly in $(C^{0,1}_0(\mathcal{D}))^*$ and $F_{\varepsilon_n}(v_n^k)\to I(\mu_k)-\zeta_{\rm max}\mu_k(\mathcal{D})$ as $n\to+\infty$. 
It then suffices to apply the diagonal argument used in the proof of Proposition \ref{limsupabscont}, Step 4,  to construct the required sequence. 
\vskip5pt

\noindent{\it Step 1.} We first consider the case where $\Sigma$ is a segment in $\DD$. Without loss of generality we may assume that $\Sigma=[a,b]\times\{0\}\subset\!\subset \DD$ 
for some $a,b\in\RR$ with $a<b$. Assume in addition that $\Sigma^\prime:={\rm supp}\,\mu\subset\!\subset ]a,b[\times\{0\}$. 
We shall regularize the measure $\mu$ using the following standard procedure.  Consider a smooth function $\varrho\in C^\infty(\RR)$ such 
that $\varrho\geq 0$, ${\rm supp}\,\varrho \subset [-1,1] $ and $\int_\RR\varrho =1$. For a positive integer $k$ and $x=(x_1,x_2)\in\mathbb{R}^2$, we introduce $\varrho_k(x):=k\varrho(kx_1)$, 
and we define 
$$g_k(x):=\int_{\Sigma^\prime} \varrho_k(x-y)\,d\mu(y) \,.$$
By construction, the function $g_k$ is nonnegative,  smooth and supported by $[a,b]\times\RR$ for $k$ large enough.  
Next we define for $k$ large the measure 
 $$\mu_k:=g_k(x)\,\mathcal{H}^1\res\Sigma\,.$$
One may easily check that $\mu_k\in H^{-1}(\DD)$, $\mu_k(\DD)\to\mu(\DD)$ and that $\mu_k\rightharpoonup\mu$ weakly* in the sense of measures on $\DD$ 
as $k\to+\infty$. In particular, $\mu_k\to\mu$ 
strongly in $(C_0^{0,1}(\mathcal{D}))^*$.

We claim that \eqref{convapproxpot} holds. Indeed, using Fubini's theorem we first derive that
\begin{multline*}
I(\mu_k)=\frac{1}{2}\iint_{\Sigma\times\Sigma}G(z,z^\prime)\,d\mu_k(z)d\mu_k(z^\prime) =\\
= \frac{1}{2}\iint_{\Sigma^\prime\times\Sigma^\prime} \bigg( \iint_{(\Sigma\cap B_{\frac{1}{k}}(x))\times(\Sigma\cap B_{\frac{1}{k}}(y))} G(z,z^\prime) \varrho_k(z-x)\varrho_k(z^\prime-y)\,
d\mathcal{H}^1_{z}d\mathcal{H}^1_{z^\prime}\bigg)\,d\mu(x)d\mu(y)\,.
\end{multline*}
Next we observe that for $k$ large enough, we have $\Sigma\cap B_{\frac{1}{k}}(x)=(x_1,0)+J_k$  
with $J_k:=(\frac{-1}{k},\frac{1}{k})\times\{0\}$ for every $x=(x_1,x_2)\in\Sigma^\prime$. Changing variables in $(z,z^\prime)$ and using Fubini's theorem again, we obtain
\begin{align*}
I(\mu_k)&=\frac{1}{2} \iint_{\Sigma^\prime\times\Sigma^\prime} \bigg( \iint_{J_k\times J_k} 
G(x+\xi,y+\xi^\prime) \varrho_k(\xi)\varrho_k(\xi^\prime)\,d\mathcal{H}^1_{\xi}d\mathcal{H}^1_{\xi^\prime}\bigg)\,d\mu(x)d\mu(y)\\
&= \frac{1}{2} \iint_{J_k\times J_k} \bigg( 
 \iint_{\Sigma^\prime\times\Sigma^\prime} G(x+\xi,y+\xi^\prime) \, d\mu(x)d\mu(y) \bigg) \varrho_k(\xi)\varrho_k(\xi^\prime)\, d\mathcal{H}^1_{\xi}d\mathcal{H}^1_{\xi^\prime}\\
 &= \frac{1}{2} \iint_{J_k\times J_k} \bigg( \iint_{\DD\times\DD} G(x,y) \, d(\tau_{\xi}\mu)(x)d(\tau_{\xi^\prime}\mu)(y) \bigg)  \varrho_k(\xi)\varrho_k(\xi^\prime)\, d\mathcal{H}^1_{\xi}d\mathcal{H}^1_{\xi^\prime}\,.
\end{align*}
From the Green representation of $h_{\tau_{\xi}\mu}$ we infer that for every $(\xi,\xi^\prime)\in J_k\times J_k$, 
$$\iint_{\DD\times\DD} G(x,y) \, d(\tau_{\xi}\mu)(x)d(\tau_{\xi^\prime}\mu)(y)= \int_{\DD}\big(\nabla h_{\tau_{\xi}\mu}\big)\cdot \big(\nabla h_{\tau_{\xi^\prime}\mu}\big) \,dx \,.$$
Then from Lemma \ref{conttrans} we deduce that  the function 
$$\Theta:(\xi,\xi^\prime)\mapsto \int_{\DD}\big(\nabla h_{\tau_{\xi}\mu}\big)\cdot \big(\nabla h_{\tau_{\xi^\prime}\mu}\big) \,dx$$
is continuous on $B_\delta(0)\times B_\delta(0)$ for some $0<\delta<\dist(\Sigma,\partial \DD)$. Therefore, 
$$\lim_{k\to+\infty}I(\mu_k)=\lim_{k\to+\infty} \frac{1}{2} \iint_{J_k\times J_k}\Theta(\xi,\xi^\prime) \varrho_k(\xi)\varrho_k(\xi^\prime)\, d\mathcal{H}^1_{\xi}d\mathcal{H}^1_{\xi^\prime}=\frac{1}{2}\Theta(0,0) =I(\mu)\,,$$
and \eqref{convapproxpot} is proved. 
\vskip5pt

\noindent{\it Step 2.} We now consider the case where $\Sigma$ is a smooth embedded arc. We further assume that there exists a $C^1$-diffeomorphism $\Phi:\DD\to \DD$ such that $\Phi(x)=x$ in a neighborhood 
of $\partial \DD$ and $\bar\Sigma:=\Phi(\Sigma)$ is a segment compactly included in $\DD$. Let $\mu$ be a nonnegative Radon measure in $H^{-1}(\DD)$ whose support is compactly included in 
the relative interior of $\Sigma$. Denote by $\bar \mu$ the push-forward of $\mu$ through $\Phi$, {\it i.e.}, $\bar\mu:=\Phi_\#\mu$. Then ${\rm supp}\,\bar\mu$ is compactly included in 
the relative interior of $\bar\Sigma$ and $\bar\mu\in H^{-1}(\DD)$. Indeed, we easily check that 
\begin{align*}
I(\bar\mu)=\frac{1}{2}\iint_{\bar\Sigma\times\bar\Sigma}G(x,y)\,d\bar\mu(x)d\bar\mu(y) &=\frac{1}{2}\iint_{\Sigma\times\Sigma}G\big(\Phi(x),\Phi(y)\big)\,d\mu(x)d\mu(y)\\
&\leq  \frac{1}{4\pi}\iint_{\Sigma\times\Sigma}\ln\big|\Phi(x)-\Phi(y)\big|\,d\mu(x)d\mu(y)+C\\
&\leq \frac{1}{4\pi}\iint_{\Sigma\times\Sigma}\ln|x-y|\,d\mu(x)d\mu(y)+C<+\infty\,,
\end{align*}
where we have used \eqref{Green} and the constant $C$ only depends on $\Sigma$, $\mu(\Sigma)$ and $\|\nabla\Phi\|_{L^\infty(\DD)}$. 

Therefore we can apply Step 1 to $\bar\mu$ to find a sequence of measures $\{\bar \mu_k\}_{k\in\NN}$ of the form \eqref{abscontmeas} such that ${\rm supp}\,\bar\mu_k\subset \bar\Sigma$, 
$\bar\mu_k(\DD)\to\bar\mu(\DD)$, $\bar\mu_k\rightharpoonup \bar\mu$ 
weakly* as measures on $\DD$, and $I(\bar\mu_k)\to I(\bar\mu)$ as $k\to+\infty$. 
Then we set $\mu_k:= (\Phi^{-1})_\#\bar\mu_k$ for every integer $k$.  Observe that $\mu_k$ is of the form \eqref{abscontmeas}. Indeed, writing 
$\bar\mu_k=\bar g_k(x)\,\mathcal{H}^1\res\bar\Sigma$ with $\bar g_k\in C^0(\bar\Sigma)$, the area formula (see {\it e.g.} \cite{AFP}) yields
$$\mu_k=\bar g_k\circ\Phi (x)|\nabla_\tau\Phi(x)|\,\mathcal{H}^1\res\Sigma\,, $$
where  $\nabla_\tau\Phi$ denotes the tangential gradient of $\Phi$ along $\Sigma$. 
Then one may check that $\mu_k(\DD)\to\mu(\DD)$, $\mu_k\rightharpoonup \mu$ weakly* as measures on $\DD$ as $k\to+\infty$. 

We claim that \eqref{convapproxpot} holds. First write 
\begin{align}
\nonumber I(\mu_k)=\,&\frac{1}{2}\iint_{\bar\Sigma\times\bar\Sigma} G\big(\Phi^{-1}(x),\Phi^{-1}(y)\big)\,d\bar\mu_k(x)d\bar\mu_k(y) \\
\nonumber =\,&\frac{1}{2}\iint_{\bar\Sigma\times\bar\Sigma} S\big(\Phi^{-1}(x),\Phi^{-1}(y)\big)\,d\bar\mu_k(x)d\bar\mu_k(y)\\
\nonumber&\,+\frac{1}{4\pi} \iint_{\bar\Sigma\times\bar\Sigma} \ln\bigg(\frac{1}{|\Phi^{-1}(x)-\Phi^{-1}(y)|}\bigg)\,d\bar\mu_k(x)d\bar\mu_k(y)\\
\label{Ik+IIk} =&: I_k+II_k\,,
\end{align}
where $S$ denotes the {\it regular part} of the Green function $G$, {\it i.e.},
$$S(x,y):=G(x,y)+\frac{1}{2\pi}\ln|x-y|$$ 
(which is a locally smooth function on $\DD\times\DD$). 
Since $\bar\mu_k$ converges weakly* as measures to $\bar \mu$, we have 
\begin{equation}\label{convprodtens}
\bar\mu_k\otimes\bar\mu_k\rightharpoonup \bar\mu\otimes \bar\mu\quad\text{weakly* as measures on $\DD\times\DD$}\,,
\end{equation}
and we deduce that  
\begin{equation}\label{limIk}
\lim_{k\to+\infty} I_k
= \frac{1}{2}\iint_{\bar\Sigma\times\bar\Sigma} S\big(\Phi^{-1}(x),\Phi^{-1}(y)\big)\,d\bar\mu(x)d\bar\mu(y)=  \frac{1}{2}\iint_{\Sigma\times\Sigma} S(x,y)\,d\mu(x)d\mu(y)\,.
\end{equation}

Next we consider a decreasing sequence  $\alpha_n\to 0$.  
For every integer $n$, we introduce a smooth cut-off function $\chi_n\in C^\infty(\overline\DD\times \overline\DD)$ such that 
$0\leq \chi_n\leq 1$, $\chi_n(x,y)=0$ if $ |\Phi^{-1}(x)-\Phi^{-1}(y)|\geq \alpha_n$, 
and $\chi_n(x,y)=1$ if $ |\Phi^{-1}(x)-\Phi^{-1}(y)|\leq \alpha_{n+1}$. 
Note since $\bar\mu\in H^{-1}(\DD)$, the measure $\bar\mu$ has no atoms, and hence 
$\bar\mu\otimes\bar\mu$ does not charge the diagonal 
$\{x=y\}\cap\DD\times\DD$. Consequently, $\chi_n\to 0$ $\bar\mu\otimes\bar\mu$--a.e. in $\DD\times \DD$. 
Then write 
\begin{multline}\label{decompIIk}
II_k=\frac{1}{4\pi}\iint_{\bar\Sigma\times\bar\Sigma}\chi_n(x,y)\ln\bigg(\frac{1}{|\Phi^{-1}(x)-\Phi^{-1}(y)|}\bigg)\,d\bar\mu_k(x)d\bar\mu_k(y)\,+\\
+\frac{1}{4\pi}\iint_{\bar\Sigma\times\bar\Sigma}(1-\chi_n(x,y))\ln\bigg(\frac{1}{|\Phi^{-1}(x)-\Phi^{-1}(y)|}\bigg)\,d\bar\mu_k(x)d\bar\mu_k(y) =:III_k^n+IV_k^n\,.
\end{multline}
By the choice of $\chi_n$ and \eqref{convprodtens}, we have  for every $n$,
$$\lim_{k\to+\infty} IV_k^n= \frac{1}{4\pi}\iint_{\bar\Sigma\times\bar\Sigma}(1-\chi_n(x,y))\ln\bigg(\frac{1}{|\Phi^{-1}(x)-\Phi^{-1}(y)|}\bigg)\,d\bar\mu(x)d\bar\mu(y)\,.$$
Next observe that 
\begin{equation}\label{equivdiffeo}
\frac{C_1}{|x-y|}\leq \frac{1}{\big|\Phi^{-1}(x)-\Phi^{-1}(y)\big|}\leq \frac{C_2}{|x-y|}\quad\text{for every $(x,y)\in\DD\times\DD$, $x\not=y$}\,,
\end{equation}
for some constants $C_1>0$ and $C_2>0$ independent of $x$ and $y$.  Since $I(\bar\mu)<+\infty$, estimate \eqref{Green} tells us that 
the function $\ln|x-y|$ belongs to $L^1(\DD\times\DD,\bar\mu\otimes\bar\mu)$. Therefore we may apply 
the dominated convergence theorem to derive  
\begin{equation}\label{limIVkn}
\lim_{n\to+\infty}\lim_{k\to+\infty} IV_k^n= \frac{1}{4\pi}\iint_{\bar\Sigma\times\bar\Sigma}\ln\bigg(\frac{1}{|\Phi^{-1}(x)-\Phi^{-1}(y)|}\bigg)\,d\bar\mu(x)d\bar\mu(y)\,.
\end{equation}

Let us now treat the term $III_k^n$. We first  deduce  from \eqref{equivdiffeo} that 
\begin{multline}\label{IIIink}
\frac{1}{4\pi}\iint_{\bar\Sigma\times\bar\Sigma}\chi_n(x,y)\ln\bigg(\frac{C_1}{|x-y|}\bigg)\,d\bar\mu_k(x)d\bar\mu_k(y)\leq  III_k^n \leq\\
\leq \frac{1}{4\pi}\iint_{\bar\Sigma\times\bar\Sigma}\chi_n(x,y)\ln\bigg(\frac{C_2}{|x-y|}\bigg)\,d\bar\mu_k(x)d\bar\mu_k(y)\,.
\end{multline}
Since the function 
\begin{multline*}
G(x,y)-\frac{\chi_n(x,y)}{2\pi}\ln\bigg(\frac{C_i}{|x-y|}\bigg)=S(x,y)+\frac{1-\chi_n(x,y)}{2\pi}\ln|x-y|-
\frac{\chi_n(x,y)}{2\pi}\ln(C_i)
\end{multline*}
is locally smooth in $\DD\times\DD$ and $I(\bar\mu_k)\to I(\bar\mu)$, we infer from \eqref{convprodtens} that  
\begin{multline}\label{intlogdiag}
\lim_{k\to +\infty} \frac{1}{4\pi}\iint_{\bar\Sigma\times\bar\Sigma}\chi_n(x,y)\ln\bigg(\frac{C_i}{|x-y|}\bigg)\,d\bar\mu_k(x)d\bar\mu_k(y)= \\
 \frac{1}{4\pi}\iint_{\bar\Sigma\times\bar\Sigma}\chi_n(x,y)\ln\bigg(\frac{C_i}{|x-y|}\bigg)\,d\bar\mu(x)d\bar\mu(y)\quad\text{for $i=1,2$}\,.
\end{multline}
Using that $\chi_n\to 0$ $\bar\mu\otimes\bar\mu$--a.e. and $\ln|x-y|$ belongs to $L^1(\DD\times\DD,\bar\mu\otimes\bar\mu)$, we infer as previously that 
\begin{equation}\label{intlogdiagvan}
\lim_{n\to+\infty} \frac{1}{4\pi}\iint_{\bar\Sigma\times\bar\Sigma}\chi_n(x,y)\ln\bigg(\frac{C_i}{|x-y|}\bigg)\,d\bar\mu(x)d\bar\mu(y)=0\quad\text{for $i=1,2$}\,.
\end{equation}
Combining \eqref{IIIink}, \eqref{intlogdiag} and \eqref{intlogdiagvan} we derive 
$$\lim_{n\to+\infty}\liminf_{k\to+\infty} III_k^n =\lim_{n\to+\infty}\limsup_{k\to+\infty} III_k^n =0 \,,$$
which yields together with \eqref{decompIIk} and \eqref{limIVkn}, 
\begin{multline}\label{limIIk}
\lim_{k\to+\infty} II_k= \frac{1}{4\pi} \iint_{\bar\Sigma\times\bar\Sigma} \ln\bigg(\frac{1}{|\Phi^{-1}(x)-\Phi^{-1}(y)|}\bigg)\,d\bar\mu(x)d\bar\mu(y)\\
=\frac{1}{4\pi} \iint_{\Sigma\times\Sigma} \ln\bigg(\frac{1}{|x-y|}\bigg)\,d\mu(x)d\mu(y)\,.
\end{multline}
Then \eqref{convapproxpot} follows gathering \eqref{Ik+IIk}, \eqref{limIk} and \eqref{limIIk}. 
\vskip5pt

\noindent{\it Step 3.} We now consider the general $\Sigma$ case. If $\Sigma$ is an embedded arc, we may assume without loss of 
generality that $\Sigma \subset \Sigma^\prime$ for some $C^2$--Jordan curve $\Sigma^\prime$ compactly included in $\DD$. 
Hence it suffices to consider the case where $\Sigma$ is a Jordan curve. We shall use the following lemma. Its proof is postponed at the end of the section. 

\begin{lemma}\label{constructdiffeo}
Assume that $\Sigma$ is a $C^2$-Jordan curve. Then there exists 
$\delta_1>0$ such that for every $x_0\in\Sigma$, 
there exists a 
$C^1$-diffeomorphism $\Phi:\DD\to\DD$ satisfying $\Phi(x)=x$ in a neighborhood of 
$\partial \DD$ and such that $\Phi(\Sigma \cap \overline B_{\delta_1}(x_0))$ 
is a segment compactly included in $\DD$. 
\end{lemma}

Now let $\gamma:[0,1]\to\Sigma$ be a constant speed parametrization of $\Sigma$. 
Let $N$ be a positive integer to be chosen and set $t_n=n/N$ for $n=0,\ldots,N$, and 
$$\Sigma_n:=\gamma([t_{n-1},t_n])\quad\text{for $n=1,\ldots,N$.}$$ 
We choose $N$ in such a way that ${\rm diam}(\Sigma_n)\leq\delta_1$ for each $n$, where 
the constant $\delta_1$ is given by Lemma~\ref{constructdiffeo}. Setting $x_n=\gamma((t_{n-1}+t_n)/2)$ for $n=1,\ldots,N$, we can apply Lemma~\ref{constructdiffeo} to each $x_n$ 
to find a $C^1$-diffeomorphism $\Phi_n:\DD\to\DD$ such that $\Phi_n(\Sigma_n)$ 
is a segment compactly included in $\DD$, and $\Phi_n(x)=x$ in a neighborhood of $\partial\DD$.

Let $\mu$ be an arbitrary nonnegative Radon measure in $H^{-1}(\DD)$ supported by $\Sigma$. 
Consider a decreasing sequence $\alpha_k\to0$ and define for $k$ large enough,
$$\Sigma^k_n:=\gamma([t_{n-1}+\alpha_k,t_n-\alpha_k])\,,\quad\mu_n^k:=\mu\res \Sigma^k_n\quad\text{for $n=1,\ldots,N$.}$$  
Oviously $\mu_n^k\in H^{-1}(\DD)$ with ${\rm supp}\,\mu_n^k\subset \Sigma^k_n$. 
Applying Step 2 for each 
$n$ and $k$, we find a sequence of measures $\{\mu^k_{n,m}\}_{m\in\NN}$ of the form \eqref{abscontmeas} such that ${\rm supp}\,\mu^k_{n,m}\subset \Sigma_n$, 
$\mu^k_{n,m}(\DD)\to\mu^k_{n}(\DD)$, $\mu^k_{n,m}\rightharpoonup \mu^k_{n}$  
weakly* as measures on $\DD$, and $I(\mu^k_{n,m})\to I(\mu^k_{n})$ as $m\to+\infty$. Define the measures 
$$ \mu^k_m:=\sum_{n=1}^N \mu^k_{n,m}\quad \text{and}\quad \mu^k:=\sum_{n=1}^N\mu^k_n=\mu\res(\cup_n \Sigma_n^k)\,.$$
Then $\mu^k_{m}(\DD)\to\mu^k(\DD)$ and $\mu^k_{m}\rightharpoonup \mu^k$  
weakly* as measures on $\DD$ as $m\to+\infty$. In addition, from Lemma \ref{convIgen} we infer that $h_{\mu^k_{n,m}}\to h_{\mu^k_n}$ strongly in $H^1(\DD)$ for 
every integers $n$ and $k$. Hence 
\begin{multline*}
I(\mu^k_m)=\frac{1}{2}\sum_{i,j=1}^N\iint_{\DD\times\DD}G(x,y)\,d\mu^k_{i,m}d\mu^k_{j,m}=\frac{1}{2}\sum_{i,j=1}^N\int_{\DD}\nabla h_{\mu^k_{i,m}}\cdot \nabla h_{\mu^k_{j,m}}\,dx\\
\mathop{\longrightarrow}\limits_{m\to+\infty} \frac{1}{2}\sum_{i,j=1}^N\int_{\DD}\nabla h_{\mu^k_{i}}\cdot \nabla h_{\mu^k_{j}}\,dx=\frac{1}{2}\sum_{i,j=1}^N\iint_{\DD\times\DD}G(x,y)\, d\mu^k_id\mu_j^k=I(\mu^k)\,.
\end{multline*}

Next recall that $\mu$ is atomless. Hence, 
by monotone convergence we have $\mu_k(\DD)\to \mu(\DD)$ and $I(\mu^k)\to I(\mu)$ as $k\to+\infty$, as well as the weak*  convergence of $\mu^k$ to $\mu$. Consequently, 
\begin{multline*}
\lim_{k\to+\infty}\lim_{m\to+\infty}|\mu^k_m(\DD)-\mu(\DD)|=\lim_{k\to+\infty}\lim_{m\to+\infty}|I(\mu^k_m)-I(\mu)|\\
=\lim_{k\to+\infty}\lim_{m\to+\infty}\|\mu^k_m-\mu\|_{(C^{0,1}_0(\DD))^*}=0
\end{multline*}
(here we use again the compact embedding $(C^0_0(\mathcal{D}))^*\hookrightarrow(C_0^{0,1}(\mathcal{D}))^*$), 
and the conclusion follows for a suitable diagonal sequence $\mu_k=\mu_{m_k}^k$.  
\end{proof}

\noindent{\it Proof of Lemma \ref{constructdiffeo}.} By assumption on $\Sigma$, 
there exists $\delta_0>0$ such that for every 
$x_0\in\Sigma$, $\Sigma \cap \overline B_{2\delta_0}(x_0)$ is the graph of a $C^2$-function and $
\overline B_{2\delta_0}(x_0)\subset \DD$.  
Now fix $x_0\in\Sigma$ and write every $x\in \DD$ as $x=x_0 +s\tau +t \tau^\perp$ where $\tau$ denotes a unit tangent vector to $\Sigma$ at $x_0$. Then 
$\Sigma \cap \overline B_{2\delta_0}(x_0)=\{x_0+ s\tau +f(s)\tau^\perp\,,\, s\in[s_{\rm min},s_{\rm max}]\}$ for some $0>s_{\rm min}\geq-2\delta_0$, $0<s_{\rm max}\leq2\delta_0$, and  
 a $C^2$-function $f:[s_{\rm min},s_{\rm max}]\to\mathbb{R}$ satisfying $f(0)=f^\prime(0)=0$. 
 Since $\Sigma$ is $C^2$, there exists a constant $\kappa>0$ which only depends on $\Sigma$ such that $|f^{\prime\prime}(s)|\leq \kappa$ for every $s\in[s_{\rm min},s_{\rm max}]$. Hence we may 
choose $\delta_0$ smaller if necessary (uniformly with respect to $x_0$) in such a way that 
$|f^\prime|\leq 1$. Then $s_{\rm min}\leq -\delta_0$, $s_{\rm max}\geq \delta_0$ and $\Sigma \cap \overline B_{\delta}(x_0)$ is still a connected arc for any $\delta\leq 2\delta_0$.

 Set $\delta_1:=\delta_0/(2+\kappa)$. We claim that  $\Sigma \cap \overline B_{\delta_1}(x_0)$ 
 satisfies the requirement.  
Indeed, we may construct a $C^1$-diffeomorphism $\Phi:\DD\to\DD$ as follows. Consider a smooth cut-off function $\chi:\RR^2\to\RR$ such that $0\leq \chi \leq 1$, $\chi(x)=1$ if $|x|\leq \delta_1$,  
$\chi(x)=0$ if $|x|\geq2\delta_0$ and $|\nabla\chi|\leq \delta_0^{-1}$.  Then we set for $x\in\DD$, 
$$\Phi(x):=x-\chi(x-x_0)f\big((x-x_0)\cdot\tau\big)\tau^\perp \,.$$
The reader may check that $\Phi$ maps $\DD$ into $\DD$, $\Phi$ is one-to-one and defines a $C^1$-diffeomorphism. Moreover $\Phi(\Sigma \cap \overline B_{\delta_1}(x_0))=
\{x_0+s\tau\,,\,-\delta_1\leq s\leq \delta_1\}$ is a segment compactly included in $\DD$. 
\prbox
\vskip10pt

\noindent{\bf Proof of Corollary \ref{asymptmin}.} 
{\it Step 1.} For any nonnegative Radon measure $\mu$ supported by $\Sigma$ we have 
\begin{equation}\label{bdinfenmeas}
I(\mu)-\zeta_{\rm max}\mu(\DD)\geq I_*(\mu(\DD))^2-\zeta_{\rm max} \mu(\DD)\,,
\end{equation} 
and equality holds if and only if $\mu=\lambda \mu_*$ for some constant $\lambda\geq 0$. We recall that  $\mu_*$ is the unique minimizer of $I$ among all probability measures supported by $\Sigma$ 
and that $I_*:=I(\mu_*)$.  
The existence and uniqueness of $\mu_*$ is classical, and we refer to \cite{ST} for further details. 
Optimizing \eqref{bdinfenmeas} with respect to $\lambda$ for measures of the form $\mu=\lambda \mu_*$, we derive that 
$\frac{\zeta_{\rm max}}{2I_*}\mu_*$ is the unique minimizer of $\mu\mapsto I(\mu)-\zeta_{\rm max}\mu(\DD)$ over all nonnegative Radon measures supported by $\Sigma$. 
\vskip5pt

\noindent{\it Step 2.} Let $\varepsilon_n\to 0^+$ be an arbitrary sequence. The existence of a minimizer $u_n$ of $F_{\varepsilon_n}$ is 
classical and follows from standard arguments based on coercivity 
and lower semicontinuity properties of $F_{\varepsilon_n}$. We first observe that $F_{\varepsilon_n}(u_n)\leq F_{\varepsilon_n}(1)=0$. 
Hence, by Theorem~\ref{mainthm}, there exists a subsequence $\{\varepsilon_{n_k}\}$ such that 
$$\frac{1}{\omega_{n_k}}j(u_{n_k})\to \mu_0\quad\text{strongly in $(C^{0,1}_0(\DD))^*$ as $k\to+\infty$,}$$
for some nonnegative Radon measure $\mu_0\in H^{-1}(\DD)$ supported by $\Sigma$. Moreover, 
\begin{equation}\label{liminfminim}
\liminf_{k\to+\infty}F_{\varepsilon_{n_k}}(u_{n_k})\geq I(\mu_0)-\zeta_{\rm max}\mu_0(\DD)\,.
\end{equation}
On the other hand, by Theorem \ref{mainthm}, any nonnegative Radon measure $\mu\in H^{-1}(\DD)$ supported by $\Sigma$ can be 
strongly approximated in $(C^{0,1}_0(\DD))^*$ by some sequence $\{\omega_{n_k}^{-1}j(v_k)\}$ 
with $\{v_k\}\subset H^1(\DD;\CC)$ satisfying 
 $$\lim_{k\to+\infty}F_{\varepsilon_{n_k}}(v_{k})= I(\mu)-\zeta_{\rm max}\mu(\DD)\,. $$
Since $F_{\varepsilon_{n_k}}(u_{n_k})\leq F_{\varepsilon_{n_k}}(v_{k})$ we infer that $\mu_0$ minimizes $\mu\mapsto I(\mu)-\zeta_{\rm max}\mu(\DD)$ 
over all nonnegative Radon measures supported by $\Sigma$. 
Consequently, $\mu_0= \frac{\zeta_{\rm max}}{2I_*}\mu_*$ and the $\liminf$ in \eqref{liminfminim} is actually a limit. Then the result along the 
full sequence $\{\varepsilon_{n}\}$ follows from a standard argument 
on the uniqueness of the limit.\prbox

\section{$\Gamma$--convergence analysis for annular domains}\label{annulus}           

In this section we briefly show how to extend the above techniques to the case of a multiply connected domain.  The method we outline here may be applied 
for any finite number of holes (see \cite{AB1}), but for simplicity we restrict to domains which are topological annuli.  
Let $\DD$ denote a simply connected domain in $\RR^2$ with smooth boundary, and $\BB\subset\!\subset\DD$ a smooth, simply connected domain compactly contained inside $\DD$.  
Then let $\AA:=\DD\setminus\BB$.  For  $u\in H^1(\AA; \CC)$ we define the functional
$$ J_\varepsilon (u) := \int_\AA \left\{ \frac12 |\nabla u|^2 + \frac{1}{4\eps^2} (|u|^2 -1)^2- \Omega_\eps V(x)\cdot j(u) \right\}\,dx\,.
$$
Here the given vector field $V:\RR^2\to\RR^2$ is assumed (for simplicity) to be locally Lipschitz continuous. 
We are interested in the asymptotic behavior of $J_\varepsilon $ as $\eps\to 0$, with an angular speed $\Omega_\eps$ as in \eqref{omega}.


\subsection{Asymptotic vorticity of the hole}

For multiply connected domains, the highest order term in an expansion of the minimal energy is partially due to the turning of the phase of a minimizer around the holes.  
The first step in studying vortices in the interior is to identify the {\it asymptotic vorticity of the hole}, and then split the energy into contributions from the hole and from the interior. 
To this purpose we first study the minimization of the functional $J_\varepsilon$ 
over $\mathbb{S}^1$--valued maps. Observe that for  $\mathbb{S}^1$--valued maps, 
the functional $J_\varepsilon$ only depends on the angular speed 
and not anymore on $\varepsilon$ itself, {\it i.e.}, 
for every  $u\in H^1(\AA;\mathbb{S}^1)$, 
$$ J_\varepsilon(u)=\mathcal{H}_\Omega (u) := \int_\AA \left\{ \frac12 |\nabla u|^2 - \Omega V(x)\cdot j(u) \right\}\,dx\,, $$
with $\Omega=\Omega_\varepsilon$. 
We are therefore interested in minimizing $\mathcal{H}_\Omega$ over the class $H^1(\AA;\mathbb{S}^1)$, and here $\Omega>0$ could be any positive parameter. 
It is well known that maps in $H^1(\AA;\mathbb{S}^1)$ are classified by their topological degree, {\it i.e.}, their winding number around the hole $\BB$. Hence, minimizing first in each 
homotopy class and then choosing the lowest energy level, one reaches the minimum of the energy 
of $\mathcal{H}_\Omega$, {\it i.e.}, 
\begin{equation}\label{mindeg}
\min \mathcal{H}_\Omega=\min_{d\in\ZZ} g(d,\Omega)\,,
\end{equation}
where 
\begin{equation}\label{defgdom}
g(d,\Omega):=\min\big\{\mathcal{H}_\Omega (u)\,:\; u\in H^1(\AA;\mathbb{S}^1)\,,\;{\rm deg}\,u =d\big\}\,. 
\end{equation}
Concerning the minimization problem \eqref{defgdom}, we have the following result. 

\begin{proposition}\label{mindfix}
For every $d\in\ZZ$, the minimization problem \eqref{defgdom} admits a unique solution $u_d$ up to a (complex) multiplicative constant of modulus one. Moreover,
\begin{equation}\label{exprgphid}
 g(d,\Omega)=\frac{1}{2}\int_{\AA}\left\{|\nabla \Phi_d|^2-\Omega^2|V|^2\right\}\,dx\,,
 \end{equation}
where $\Phi_d$ is the unique solution of the linear equation
\begin{equation}\label{constrphase}
\begin{cases}
-\Delta \Phi_d =\Omega\, \curl V & \text{in } \AA\,,\\
\Phi_d = 0  & \text{on }\partial \DD\,,\\
\Phi_d = {\rm const.} & \text{on }\partial \BB\,,  \\[8pt]
\displaystyle \int_{\partial \BB} \frac{\partial  \Phi_d }{\partial \nu} = 2\pi d -\Omega \int_{\partial \BB} V\cdot\tau
\end{cases} 
\end{equation}
\end{proposition}

\begin{proof}  We follow here some of the arguments in \cite[Chap. 1]{BBH2}, and we provide some details for the reader conveniance. 
\vskip5pt

\noindent{\it Step 1.} We  claim that for any $u\in H^1(\AA;\mathbb{S}^1)$ such that ${\rm deg}\,u =d$, we have 
$$\mathcal{H}_\Omega (u)\geq \frac{1}{2}\int_{\AA}\left\{|\nabla \Phi_d|^2-\Omega^2|V|^2\right\}\,dx\,.$$
Indeed, we first observe that $\curl j(u)=0$ since $u$ is $\mathbb{S}^1$-valued. On the other hand, 
$\nabla^\perp \Phi_d+\Omega V $ is also curl-free   
and 
$$ \int_{\partial \BB} \big(j(u) -\nabla^\perp \Phi_d-\Omega V\big)\cdot \tau = 2\pi d- \int_{\partial \BB}  \frac{\partial \Phi_d}{\partial \nu} -\Omega \int_{\partial \BB} V\cdot\tau = 0\,,$$
so that we can find a scalar function $H\in H^1({\cal A})$ such that $j(u)= \nabla H +\nabla^\perp \Phi_d+\Omega V$. 
Since $u$ is $\mathbb{S}^1$-valued, we have $|\nabla u|^2=|j(u)|^2$, and thus  
\begin{align*}
\mathcal{H}_\Omega (u) &=\frac{1}{2}\int_{\AA}\left\{|j(u)-\Omega V|^2-\Omega^2|V|^2\right\}\,dx\\
& = \frac{1}{2}\int_{\AA}\left\{|\nabla \Phi_d|^2-\Omega^2|V|^2\right\}\,dx+\frac{1}{2}\int_{\AA}|\nabla H|^2dx + \int_\AA \nabla^\perp \Phi_d\cdot \nabla H\,dx\,.
\end{align*}
Then, using the fact that the function $\Phi_d$ is constant on $\partial \AA$, an integration by parts yields $\int_\AA \nabla^\perp \Phi_d\cdot \nabla H\,dx=0$  and the claim follows. 
\vskip5pt
	
\noindent{\it Step 2.} 	 We claim that there exists $u_d\in H^1(\AA;\mathbb{S}^1)$ such that ${\rm deg}\,u_d =d$ and 
$$j(u_d)= \nabla^\perp \Phi_d+\Omega V\,.$$
Indeed, since ${\rm curl}(\nabla^\perp \Phi_d+\Omega V)=0$ and 
$$\frac{1}{2\pi} \int_{\partial \BB}   (\nabla^\perp \Phi_d +\Omega V)\cdot\tau \,=d\in \mathbb{Z}\,,$$
we may locally define a scalar function $\psi$ in $\AA$ such that 
$$\nabla\psi = \nabla^\perp \Phi_d+\Omega V\,.$$
Then $u_d:={\rm exp}(i\psi)$ is well defined and satisfies the required properties. Clearly the construction of $u_d$ is unique modulo a constant phase, and the proof is complete.  
\end{proof}

In order to solve problem \eqref{mindeg}, it now suffices to express \eqref{exprgphid} explicitely in 
terms of the integer $d$.  To this purpose, we first introduce the solution $\xi$ of the linear problem   
\begin{equation}\label{fctcap}
\begin{cases}
\Delta\xi=0 & \text{in }\AA\,,\\
\xi =0 & \text{on }\partial \DD\,,\\ 
\xi=1 & \text{on }\partial \BB\,.
\end{cases}
\end{equation}
The function $\xi$ is smooth in $\overline\AA$ and $0\leq \xi\leq 1$ by the maximum principle. 
Moreover, the Dirichlet energy of $\xi$ is  the so-called $H^1$-capacity of $\BB$ inside $\DD$ which we denote by ${\rm cap}(\BB)$, {\it i.e.}, 
\begin{equation}\label{capa}
{\rm cap}(\BB):=\int_\AA|\nabla\xi|^2\,dx=- \int_{\partial \BB}\frac{\partial \xi}{\partial \nu}\,>0\,.
\end{equation}
Next we consider the unique  solution $\zeta$ of
\begin{equation}\label{fctpot}
\begin{cases}
-\Delta \zeta = \curl V & \text{in } \AA\,,\\
\zeta = 0  & \text{on }\partial \AA\,.
\end{cases} 
\end{equation}
From the Lipschitz assumption on $V$ and standard elliptic regularity, we infer that $\zeta$ belongs to $C_0^{1,\alpha}(\overline\AA)$ for every $0\leq \alpha <1$. We set 
$$\gamma_V:=\int_{\partial\DD}\bigg\{ \frac{\partial \zeta}{\partial \nu} 
+V\cdot\tau\bigg\}\,.$$ Observing that \eqref{fctpot} implies 
$$\int_{\partial\BB} \frac{\partial \zeta}{\partial \nu} =\gamma_V-\int_{\partial \BB} V\cdot \tau\,,$$ 
we find that for every integer $d$,  the function $\Phi_d$ determined by \eqref{constrphase} is explicitly given by 
\begin{equation}\label{explphi}
\Phi_d = \bigg(\frac{\gamma_V\Omega-2\pi d}{{\rm cap}(\BB)}\bigg) \xi +\Omega\, \zeta\,.
\end{equation}
Moreover, using \eqref{fctcap}, \eqref{capa} and \eqref{fctpot} we readily obtain that for every $d\in\mathbb{Z}$, 
$$\frac{1}{2}\int_{\AA}\left\{|\nabla \Phi_d|^2-\Omega^2|V|^2\right\}\,dx=\frac{|\gamma_V\Omega-2\pi d|^2}{2{\rm cap}(\BB)}-\frac{\Omega^2}{2}\int_{\AA}\left\{|V|^2-|\nabla \zeta|^2\right\}\,dx \,.$$
As a consequence, an integer $d_\Omega$ is a minimizer in \eqref{mindeg} if and only if $d_\Omega$ minimizes the function $d\in\mathbb{Z}\mapsto |\gamma_V\Omega-2\pi d|$. 
\vskip5pt

We may now state our result concerning problem \eqref{mindeg}. 

\begin{theorem}\label{thmvorthole}
Up to multiplicative constants of modulus one, the minimization problem \eqref{mindeg} admits exactly two  solutions (of distinct topological degree) 
if $\gamma_V\Omega/\pi$ is an odd integer, and a unique solution otherwise. 
Moreover, if $d_\Omega\in\ZZ$ is a minimizer in  \eqref{mindeg}, then $d_\Omega\in\big\{[\frac{\gamma_V\Omega}{2\pi}],[\frac{\gamma_V\Omega}{2\pi}]+1\big\}$ where $[\cdot]$ denotes the integer part, and 
\begin{equation}\label{minengH}
\min \mathcal{H}_\Omega = - \frac{\Omega^2}{2}\int_{\AA}\left\{|V|^2-|\nabla \zeta|^2\right\}\,dx +O(1) \quad\text{as $\Omega\to+\infty$}\,.
\end{equation}
\end{theorem}


\subsection{The $\Gamma$--convergence result}

To state the parallel $\Gamma$-convergence result for the anular domain case we must give more specific hypotheses on the potential $V$ and the angular speed $\Omega_\varepsilon$.  
In addition to the Lipschitz regularity, we assume in the sequel that $V$ satisfies the following assumptions:
\begin{itemize}
\item[(H1')]\label{H1hole} the solution $\zeta$ of \eqref{fctpot} is such that $\zeta_{\max}:=\max_{x\in\overline{\AA}}\, \zeta(x)=\max_{x\in\overline{\AA}}\, |\zeta(x)|>0$;
\vskip7pt

\item[(H2')]\label{Sighole} the set $\Sigma:=\{x\in\AA: \ \zeta(x)=\zeta_{\max}\}\subset\!\subset \AA$ is
a Jordan curve or a simple embedded arc of class $C^2$.
\end{itemize}
We note that in (H1'), the assumption that $\zeta_{\rm max}$ is achieved at positive values of $\zeta$ is not restrictive. Indeed, considering the complex conjugate of an admissible function  
replaces $V$ by $-V$ in the energy and hence $\zeta$ by $-\zeta$. 
\vskip5pt

As for the simply connected domain case, we assume that $\Omega_\varepsilon$ is near  the critical value needed for the presence of vortices which again reads 
\begin{equation}\label{omegahole}
\Om_\varepsilon=\frac{\lep}{2\zeta_{\rm max}} +  \om(\varepsilon)\,,
\end{equation}
for some positive  function $\omega$ satisfying $\omega(\varepsilon)\to+\infty$ with $\omega(\varepsilon)\leq o(|\ln\varepsilon|)$ as $\varepsilon\to 0^+$, exactly as in \eqref{omega}. 
\vskip5pt

In the sequel, for an arbitrary sequence $\varepsilon_n\to 0^+$, we will denote by $u^\star_n$ a minimizer of 
$\mathcal{H}_{\Omega_n}$, {\it i.e.}, a solution of \eqref{mindeg}, and its corresponding topological degree will be denoted by 
$d_n$. For brievety we shall also write  \eqref{explphi} as 
$$\Phi_n:=\Phi_{d_n} =\alpha_n \xi +\Omega_n \zeta\quad \text{with}\quad \alpha_n:=\frac{\gamma_V\Omega_n-2\pi d_n}{{\rm cap}(\BB)}\,.$$
We emphasize that $\alpha_n=O(1)$ as $n\to+\infty$ thanks to Theorem \ref{thmvorthole}. 
\vskip5pt

For $v\in H^1(\AA;\CC)$, we now define  
$$\overline{F}_{\varepsilon_n}(v):=\omega_n^{-2} \int_{\AA}\bigg\{\frac{|\nabla v|^2}{2}+\frac{(1-|v|^2)^2}{4\varepsilon^2_n}+\nabla^\perp\Phi_n\cdot j(v)\bigg\}\,dx\,.$$
The following proposition shows that the functional $\overline{F}_{\varepsilon_n}(\bar u^\star_n u)$ captures the energy induced by interior vorticity of a given configuration $u$.    

\begin{proposition}[Energy decomposition]\label{energdecomp}
Assume that \eqref{omegahole} holds. Let $\varepsilon_n\to 0^+$ be an arbitrary sequence and $\{u_n\}_{n\in\mathbb{N}}\subset H^1(\AA;\CC)$. Then, setting 
$v_n=\bar u^\star_nu_n\in H^1(\AA;\CC)$,  we have 
\begin{equation}\label{equivenergy}
\sup_n\, \Omega_n^{-2} J_{\varepsilon_n}(u_n)<+\infty \quad\text{if and only if}\quad \sup_n \,\Omega_n^{-2}\omega_n^{2}\overline{F}_{\varepsilon_n}(v_n)<+\infty\,.
\end{equation}
Moreover, if one of the conditions in \eqref{equivenergy} holds, then 
\begin{equation}\label{relGF}
J_{\varepsilon_n}(u_n)= \min \mathcal{H}_{\Omega_n} + \omega_n^{2}\overline{F}_{\varepsilon_n}(v_n) + o(1)\quad\text{as $n\to+\infty$}\,.
\end{equation}
\end{proposition}
\begin{proof}
Straightforward computations yield
$$\frac{|\nabla u_n|^2}{2}=|u_n|^2\frac{|\nabla u^\star_n|^2}{2} +\frac{|\nabla v_n|^2}{2} + j(u^\star_n)\cdot j(v_n)\,,$$
and 
\begin{equation}\label{reljac}
j(u_n)= |u_n|^2j(u^\star_n)+ j(v_n)\,.
\end{equation}
By the proof of Proposition \ref{mindfix}, we have $j(u^\star_n)=\nabla^\perp\Phi_n +\Omega_n V$. Hence, 
 $$J_{\varepsilon_n}(u_n)= \min \mathcal{H}_{\Omega_n} + \omega_n^{2}\overline{F}_{\varepsilon_n}(v_n) + \int_{\AA}(|u_n|^2-1)\bigg\{\frac{|\nabla u^\star_n|^2}{2} -\Omega_n V\cdot j(u^\star_n)\bigg\}\,dx\,.$$
 Then we observe that
 $$\|\nabla u^\star_n\|_{L^\infty(\AA)}= \|j(u^\star_n)\|_{L^\infty(\AA)}=\|\nabla^\perp\Phi_n +\Omega_n V\|_{L^\infty(\AA)}=O(\Omega_n)\,.$$
 Therefore, 
$$\bigg| \int_{\AA}(|u_n|^2-1)\bigg\{\frac{|\nabla u^\star_n|^2}{2} -\Omega_n V\cdot j(u^\star_n)\bigg\}\,dx\bigg| \leq 
O\left(\varepsilon_n\Omega^2_n\sqrt{E_{\varepsilon_n}(u_n,\AA)}\right)\,.$$
Since $|u_n|=|v_n|$, we also have the same estimate as above with  $E_{\varepsilon_n}(v_n,\AA)$ instead of $E_{\varepsilon_n}(u_n,\AA)$.  
Assuming that one of the conditions in \eqref{equivenergy} holds and arguing as in \eqref{contrGL}, we derive that either $\sqrt{E_{\varepsilon_n}(u_n,\AA)}\leq O(\Omega_n)$, or  
$\sqrt{E_{\varepsilon_n}(v_n,\AA)}\leq O(\Omega_n)$. Consequently, if one  of the conditions in \eqref{equivenergy} is satisfied, \eqref{relGF} holds and  the conclusion 
follows combining \eqref{relGF} with \eqref{minengH}. 
\end{proof}

For a nonnegative Radon measure $\mu$ on $\mathcal{A}$, we define 
$$ \bar{I}(\mu):=\frac{1}{2}\iint_{\AA\times \AA}
\bar{G}(x,y) \, d\mu(x)\, d\mu(y) \,, $$
where the function $\bar{G}$ denotes the Dirichlet Green's function of the domain $\mathcal{A}$, {\it i.e.}, for every $y\in\mathcal{A}$, $\bar{G}(\cdot,y)$ is the solution of 
\begin{equation}\label{defGreenhole}
\begin{cases}
-\Delta \bar{G}(\cdot,y)=\delta_y &\text{in $\mathscr{D}'(\mathcal{A})$}\,,\\
\bar{G}(\cdot,y) =0 & \text{on $\partial\mathcal{A}$}\,.
\end{cases}
\end{equation}
We may now state the $\Gamma$-convergence result for annular domains which involves 
the family of {\it ``reduced"} functionals $\{\overline{F}_\varepsilon\}_{\varepsilon>0}$. 
      
\begin{theorem}\label{thmhole}
Assume that \emph{(H1')}, \emph{(H2')} and \eqref{omegahole} hold. 
Let $\varepsilon_n\to 0^+$ be an arbitrary sequence. Then,
\begin{itemize}
\item[(i)]  for any $\{v_n\}_{n\in\NN}\subset H^1(\mathcal{D};\CC)$ satisfying $\sup_n \overline{F}_{\varepsilon_n}(v_n)<+\infty$, there exist a subsequence (not relabelled) and a 
nonnegative Radon measure $\mu$ in $H^{-1}(\mathcal{A})$ supported by $\Sigma$ such that 
\begin{equation}\label{convjachole}
\frac{1}{\omega(\varepsilon_n)}\, \curl\,j(v_n)\mathop{\longrightarrow}\limits_{n\to+\infty} \mu \quad\text{strongly in $(C^{0,1}_0(\mathcal{A}))^*$}\,;
\end{equation}
\item[(ii)] for any $\{v_n\}_{n\in\NN}\subset H^1(\mathcal{D};\CC)$ such that \eqref{convjachole} holds for some nonnegative Radon measure $\mu$ in 
$H^{-1}(\mathcal{A})$ supported by $\Sigma$, we have 
$$\liminf_{n\to+\infty}\,\overline{F}_{\varepsilon_n}(v_n)\geq \bar I(\mu)-\zeta_{\rm max}\,\mu(\mathcal{A})\,;$$
\item[(iii)] for any  nonnegative Radon measure $\mu$ in $H^{-1}(\mathcal{A})$ supported by $\Sigma$, there exists a sequence $\{v_n\}_{n\in\NN}\subset H^1(\mathcal{A};\CC)$ such that 
\eqref{convjachole} holds and 
$$\lim_{n\to+\infty} \overline{F}_{\varepsilon_n}(v_n)= \bar I(\mu)-\zeta_{\rm max}\,\mu(\mathcal{A})\,.$$
\end{itemize}
\end{theorem}

As in the simply connected case this $\Gamma$-convergence result could lead to the asymptotic description of the vorticity in $\overline{F}_\varepsilon$-global minimizers. 
Actually Theorem \ref{thmhole} combined with Proposition \ref{energdecomp} also gives the asymptotic behavior of vorticity in $J_\varepsilon$-global minimizers. The key 
observation here is that minimizers for  $J_\varepsilon$ yield  quasi-minimizers for  $\overline{F}_\varepsilon$, and conversely.  

\begin{corollary}\label{asymptminhole}
Assume that \emph{(H1')}, \emph{(H2')} and \eqref{omegahole} hold. 
Let $\varepsilon_n\to 0^+$ be an arbitrary sequence. For every integer $n\in\mathbb{N}$, let $u_n\in H^1(\mathcal{A};\CC)$ be a minimizer of $J_{\varepsilon_n}(\cdot)$. 
Then, 
$$\frac{1}{\omega(\varepsilon_n)}\, \curl\,j(u_n)\mathop{\longrightarrow}\limits_{n\to+\infty} \frac{\zeta_{\rm max}}{2 \bar I_*}\,\bar \mu_* \quad\text{strongly in $(C^{0,1}_0(\mathcal{A}))^*$}\,,$$
where $\bar \mu_*$ is the unique minimizer of $\bar I(\cdot)$ over all  probability measures supported by $\Sigma$, and $\bar I_*:=\bar I(\bar \mu_*)$. 
In addition,
\begin{equation}\label{minenergexphole}
J_{\varepsilon_n}(u_n)=-\frac{\Omega^2_n}{2}\int_\AA\big\{|V|^2-|\nabla\zeta|^2\big\}\,dx -\frac{\zeta^2_{\max}}{4\bar I_*} \omega^2_n+o(\omega_n^2)\,.
\end{equation}
\end{corollary}

We conclude this subsection with an elementary example motivated by \cite{AAB,AB1,AB2}.

\begin{example}\label{exhole}
Assume that $\DD=B_1(0)$, $\BB=B_\rho(0)$ for some $0<\rho<R$ and $V(x)=x^\perp$. Then the solution $\zeta$ of  \eqref{fctpot} is given by 
$$\zeta(x)=-\frac{|x|^2}{2}+\frac{R^2-\rho^2}{2\ln(R/\rho)} \ln|x|+\frac{\rho^2\ln R -R^2\ln\rho}{2\ln(R/\rho)}\,.$$
In particular, the set $\Sigma$ is given by the concentric circle $B_{r_*}(0)$ with 
$$r_*=\sqrt{\frac{R^2-\rho^2}{2\ln(R/\rho)} }\in(\rho,R)\,.$$
Here again, the uniform measure $\mu_*=(2\pi r_*)^{-1}d\mathcal{H}^1\res\Sigma$ turns out to be the Green equilibrium measure for $\Sigma$ in $\AA$, {\it i.e.,} $\bar I(\mu_*)=\bar I_*$. 
Indeed, one may easily check that the function 
$$\displaystyle h_*(x)=\begin{cases} 
\displaystyle\frac{\ln(R/r_*)}{2\pi(\ln(R/r_*)+\ln(r_*/\rho))} \ln(|x|/\rho) & \text{if }\rho \leq |x|\leq r_*\,,  \\[8pt]
\displaystyle \frac{\ln(r_*/\rho)}{2\pi(\ln(R/r_*)+\ln(r_*/\rho))} \ln(R/|x|) & \text{if } r_*\leq |x|\leq R\,,\\
\end{cases}
$$ 
solves $-\Delta h_*=\mu_*$ in $\AA$ with ${h_*}_{|\partial \AA}=0$. Hence $h_*(x)=\int_{\AA}\bar G(x,y)\,d\mu_*(y)$, and since $h_*$ is constant on $\Sigma$ the 
conclusion follows from Theorem II.5.12 in \cite{ST}. 
\end{example}


\subsection{Compactness of normalised weak Jacobians}\label{compacthole}

In this subsection we shall be concerned with the proof of claim {\it (i)} in Theorem \ref{thmhole}. 
We consider an arbitrary sequence $\varepsilon_n\to 0^+$. For any $\{v_n\}_{n\in\NN}\subset H^1(\mathcal{D};\CC)$ satisfying 
$\sup_n \overline{F}_{\varepsilon_n}(v_n)<+\infty$, we first derive exactly as in Lemma \ref{lem1} the estimates 
\begin{equation}\label{bdglhole}
E_{\varepsilon_n}(v_n,\AA)\leq O(|\ln \varepsilon_n|^2) \quad\text{and}\quad \|v_n\|_{L^4(\AA)}\leq O(1)\,.
\end{equation}
Hence, assuming in addition that $\{v_n\}_{n\in\NN}\subset C^1(\AA)$, we can 
 apply the vortex ball construction in Proposition~\ref{consvortexball} with $\AA$ in place of $\DD$, and $\AA_\varepsilon:=\{x\in\AA\,:\,{\rm dist}(x,\partial \AA)>\varepsilon\}$ in place of $\DD_\varepsilon$.  
We choose again $r=r_n:=|\ln\eps_n|^{-4}$, thus obtaining a finite collection of disjoint closed balls $\{\overline B(a_i^n,\rho_{i,n})\}_{i\in I_n}$ (written $B_i^n:=\overline B(a_i^n,\rho_{i,n})$), 
with associated degrees $d_{i,n}$ and total approximate vorticity  
$$D_n:=\sum_{i\in I_n}|d_{i,n}|\,, $$
as in Proposition~\ref{consvortexball}.

The first difference with the simply connected case arises in an estimate  analogue to \eqref{lem2}, 
as we must take into account an additional contribution to the potential term due to the boundary $\partial\BB$.
Since one or more of the vortex balls $B_i^n$ may intersect $\partial \AA_\varepsilon\setminus \partial\DD_\varepsilon$, we will need to perturb this boundary slightly.  
When calculating the boundary term it will be convenient to choose a level set of $\xi\,$.  For $0<s<t<1$, denote by 
\begin{equation}\label{sigma}  
\sigma_t:=\{x\in\AA: \xi(x)=t\}, \text{ \ and \ } 
   A_{s,t}:=\left\{x\in\AA: \ s<\xi<t\right\}\,.  
\end{equation}
As an easy consequence of the Maximum Principle and the Hopf boundary lemma, each curve $\sigma_t$ is smooth and 
the family $\{\sigma_t\}_{0<t<1}$ realizes a foliation of $\AA$. Then,  
for every $t\in(0,1)$  the curve $\sigma_t$ is diffeomorphic to
$\partial\BB$, and the set $A_{t,1}$ is a neighborhood
 of  $\partial\BB$ in $\AA$. 
 Now we shall choose an appropriate level set of $\xi$. Define
\begin{equation}\label{goodlevcurv}
  J_n=\big\{ t\in (0, 1)\,: \, \sigma_t\cap \left(\cup_{i\in I_n} B_i^n\right)=\emptyset\big\}\,.
\end{equation}
We note that the measure of the complement $(0,1)\setminus J_n$ is of the same order as $r_n=|\ln\eps_n|^{-4}$.  
Hence we can find $t_n\in J_n$ such that the level curve 
$\gamma_n := \sigma_{t_n}$ satisfies 
$$\varepsilon_n<\dist(\gamma_n,\partial\BB)\leq O(|\ln\eps_n|^{-3})\,,$$
and consequently $t_n=1+O(|\ln\eps_n|^{-3})$.  

This construction allows us to define to topological degree of $v_n$ around $\gamma_n$ since $v_n$ does not vanish on $\gamma_n$, {\it i.e.}, 
$$\delta_n:= {\rm deg}\bigg(\frac{v_n}{|v_n|},\gamma_n\bigg)\,.$$
Then an approximate total vorticity in $\AA$ of the configuration $v_n$ is given by $|\delta_n|+D_n$. 
\vskip3pt

We may now state the following proposition which parallels Proposition \ref{degreebound}.

\begin{proposition}\label{degreeboundhole}
Assume that \emph{(H1')} and \eqref{omegahole} hold. Let $\{v_n\}_{n\in\NN}\subset H^1(\mathcal{D};\CC)\cap C^1(\AA)$, $D_n$ 
and $\delta_n$ as above. Then  $|\delta_n|+D_{n}\leq O(\omega_n)$. 
\end{proposition}

\begin{proof}
{\it Step 1.} Arguing exactly as in the proof of \eqref{lem2}, we first derive
\begin{equation}\label{lem2hole}
\Omega_n\int_\AA \nabla^\perp \zeta \cdot j(v_n) = -2\pi\Omega_n\sum_{i\in I_n} d_{i,n}\zeta(a_i^n)  +o(1)\,.
\end{equation}
Then following essentially the proof of Lemma 3.4 in \cite{AB1}, we obtain 
\begin{equation}\label{actxi}
\int_{\AA}\nabla^\perp\xi \cdot j(v_n) = -2\pi\sum_{i\in I_n} d_{i,n} \xi(a_i^n) -2\pi t_n \delta_n +o(1)\,.
\end{equation}
Here the fact that $\gamma_n$ is a level set of $\xi$ is essential in obtaining the degree $\delta_n$ from the boundary term when integrating by parts. 
\vskip5pt

\noindent{\it Step 2.} As in the proof of Proposition \ref{degreebound},  we may assume that $\omega_n\leq O(D_n+|\delta_n|)$. From \eqref{bdglhole} and claim {\it (iv)} in
 Proposition~\ref{consvortexball}, we have $D_n\leq O(|\ln\varepsilon_n|)$. 
Then we infer from \eqref{balls3}, \eqref{lem2hole} and \eqref{actxi} that 
\begin{multline}\label{lb1hole}
O(\omega^2_n) \geq  \omega^2_n \overline{F}_{\varepsilon_n}(v_n)
\geq   \pi D_{n} \big(|\ln\varepsilon_n|-C\ln|\ln\varepsilon_n|\big)
     - 2\pi\Omega_{n}\sum_{i\in I_{n}} d_{i,n}\zeta(a^n_i) \,-2\pi |\alpha_n| D_n \\ 
-2\pi|\alpha_n||\delta_n| +  \int_{\AA\setminus\cup_{i\in I_n} B^n_i} |\nabla v_n|^2\, dx +o(1)\,,
\end{multline}
where we have used the fact that $0\leq \xi\leq 1$ and $0<t_n<1$.  

Next we consider a sequence $\eta_n\to 0$ as in \eqref{eta}, and we group the vortex balls into different classes as in the proof of Proposition \ref{degreebound}
(we refer to it for the notation). Exactly as in \eqref{lb*} and \eqref{lb+}, we derive that 
\begin{equation}\label{lb*hole}
\pi D_n^*|\ln\varepsilon_n| - 2\pi \Omega_n \sum_{i\in I_n^*} d_{i,n}\zeta(a^n_i)
  \geq  -2\pi \omega_n \zeta_{\rm max} D_n^*\,,
\end{equation}
and
\begin{equation}\label{lb+hole}
\pi D_n^+|\ln\varepsilon_n| - 2\pi \Omega_n \sum_{i\in I_n^+} d_{i,n}\zeta(a^n_i) 
\geq C\Omega_n\eta_n D_n^+\,.
\end{equation}
For negative degrees, we observe that  $\{\zeta\leq 0\}\cap \Sigma=\emptyset$ since $\zeta_{\max}=|\zeta|_{\max}>0$. Hence 
we can estimate as for the class $I_n^+$, 
\begin{multline}\label{lb-hole}
 \pi D_n^- |\ln\varepsilon_n|
     - 2\pi\Omega_n\sum_{i\in I_n^-} d_{i,n} \zeta(a^n_i)
\geq 
 \pi D_n^- |\ln\varepsilon_n|
     - 2\pi\Omega_n\sum_{i\in I_n^-,\zeta(a^n_i)\leq 0} d_{i,n} \zeta(a^n_i)\\
     \geq C\Omega_n\eta_n D_n^-\,.
\end{multline}
Inserting \eqref{lb*hole}, \eqref{lb+hole} and \eqref{lb-hole} in \eqref{lb1hole} yields
\begin{multline} \label{lb2hole}  
O(\omega_n^2) \geq 
 -\pi C_1 D_n\ln|\ln\varepsilon_n| -2\pi \omega_n \zeta_{\rm max} D_n +C\eta_n\Omega_n (D_n^++D_n^-) \\
-2\pi|\alpha_n||\delta_n|  + \frac{1}{2}\int_{\AA\setminus\cup B_i^n} |\nabla v_n|^2 dx\,. 
\end{multline}
Using the fact that $|\alpha_n|=O(1)$, we easily deduce the estimate
\begin{equation}
\label{prebdD+hole} D_n^++D_n^- \leq C \,\frac{\max\{\omega_n,\ln|\ln\varepsilon_n|\}(D_n+|\delta_n|)}{ \eta_n|\ln\varepsilon_n|}  \,, 
\end{equation}
for a constant $C>0$ independent of $n$.

We claim that 
\begin{equation}\label{estideghole}
\int_{\AA\setminus\cup B_i^n} |\nabla v_n|^2 dx \geq C(D_n+|\delta_n|)^2\,.
\end{equation}
Accepting \eqref{estideghole}, we infer from \eqref{lb2hole} that
$$(D_n+|\delta_n|)^2 -C\max\{\omega_n,\ln|\ln\varepsilon_n|\}(D_n+|\delta_n|)\leq O(\omega_n^2)\,, $$
which clearly implies
\begin{equation}\label{predegbdhole}
D_n+|\delta_n|\leq O(\max\{\omega_n,\ln|\ln\varepsilon_n|\})\,.
\end{equation}
\vskip3pt

To prove \eqref{estideghole} we introduce 
$$t_0:=1/2 \min\{t\in(0,1): \sigma_t\cap\Sigma\not=\emptyset\}\,, \quad t_1:=1/2(1+\max\{t\in(0,1): \sigma_t\cap\Sigma\not=\emptyset\})\,,$$
and 
$$  J^1_n:=\big\{ t\in (\varepsilon_n,t_0)\,: \, \sigma_t\cap \left(\cup_{i\in I_n} B_i^n\right)=\emptyset\big\}\,,\quad
J^2_n:=\big\{ t\in (t_1,t_n)\,: \, \sigma_t\cap \left(\cup_{i\in I_n} B_i^n\right)=\emptyset\big\}\,. $$
Then $\dist(\sigma_t,\Sigma)\geq C>0$ for every $t\in J_n^1\cup J_n^2$ and $\min(\mathcal{L}^1(J_n^1),\mathcal{L}^1(J_n^2))\geq C>0$ for 
a constant $C$ independent of $n$.  Next we consider
$$D_n(t):={\rm deg}\bigg(\frac{v_n}{|v_n|},\sigma_t\bigg)\,. $$
If $2|\delta_n|\geq D_n$ and  $n$ large enough, we estimate using \eqref{prebdD+hole}, 
\begin{equation}\label{estideglargedelt}
|D_n(t)|\geq \big||\delta_n|-D^+_n-D_n^-\big|\geq \frac{1}{2}|\delta_n|\geq C(D_n+|\delta_n|)\quad\text{for every $t\in J_n^2$}\,.
\end{equation}
In the opposite case $2|\delta_n|< D_n$ (and  $n$ large), we have in view of \eqref{prebdD+hole},
\begin{equation}\label{estidegsmalldelt}
 |D_n(t)|\geq \big|D_n-D_n^+-2D_n^- -|\delta_n|\big|\geq \frac{1}{4}D_n \geq C(D_n+|\delta_n|)\quad\text{for every $t\in J_n^1$}\,.
 \end{equation}
Set $\tilde v_n := v_n/|v_n|$. Using claim {\it (ii)} in Proposition~\ref{consvortexball}, the Coarea Formula and Jensen Inequality, we derive from \eqref{estideglargedelt} and \eqref{estidegsmalldelt} that 
\begin{multline*}
\int_{\AA\setminus\cup B_i^n} |\nabla v_n|^2 \,dx \geq C\int_{\AA\setminus\cup B_i^n} |\nabla \tilde v_n|^2|\nabla\xi|\, dx\geq C\int_{J_n^1\cup J_n^2}\bigg(\int_{\sigma_t}|\nabla \tilde v_n|^2\bigg) dt\\ 
\geq C \int_{J_n^1\cup J_n^2}\bigg(\int_{\sigma_t}|\tilde v_n\wedge\nabla_{\tau} \tilde v_n|^2d\mathcal{H}^1\bigg) dt\geq C \int_{J_n^1\cup J_n^2} \frac{|D_n(t)|^2}{\mathcal{H}^1(\sigma_t)}\,dt \geq C(D_n+|\delta_n|)^2\,,
\end{multline*}
and \eqref{estideghole} is proved (here we have also used the fact that $|\nabla \xi|$ does not vanish in $\overline\AA$). 
\vskip5pt

\noindent{\it Step 3.} If $\ln|\ln\varepsilon_n|=o(\omega_n)$, the conclusion follows from \eqref{prebdD+hole}  and \eqref{predegbdhole}. 
If $\omega_n=O(\ln|\ln\varepsilon_n|)$, we refine the lower 
bound by  growing the vortex balls as already performed in Step 3 of the proof of Proposition \ref{degreebound}. Following the same arguments 
with minor modifications yields the announced result, so we  omit the details.  
\end{proof}

\noindent{\bf Proof of Theorem \ref{thmhole}, claim {\it (i)}.} In view of Proposition \ref{degreeboundhole} and estimate \eqref{prebdD+hole}, 
we can follow the proof of Theorem \ref{compactness}, considering first  a sequence $\{v_n\}_{n\in\NN}\subset H^1(\mathcal{D};\CC)\cap C^1(\AA)$ and then the general case.
\prbox


\subsection{The lower bound inequality}

\noindent{\bf Proof of Theorem \ref{thmhole}, claim{\it  (ii)}.} 
We shall use the notations of Subsection \ref{compacthole}. 
Without loss of generality, we may assume that 
\begin{equation}\label{liminflimhole}
\liminf_{n\to+\infty} \overline{F}_{\varepsilon_n}(v_n)=\lim_{n\to+\infty} \overline{F}_{\varepsilon_n}(v_n)<+\infty\,.
\end{equation}
As in the proof of Theorem \ref{lwbd}, we may also assume that 
$\{u_n\}_{n\in\NN}\subset H^1(\mathcal{D};\mathbb{C})\cap C^1(\mathcal{D})$. We shall only consider the case $\ln|\ln\varepsilon_n|=o(\omega_n)$ since the other case can be 
completed as we already pursued in the proof of Theorem \ref{lwbd}. 

By Proposition \ref{degreeboundhole} and \eqref{actxi}, we have 
$$\bigg| \alpha_n \int_{\AA}\nabla^\perp\xi \cdot j(v_n) \,dx\bigg| =O(\omega_n)\,.$$
Then we can argue exactly as in \eqref{saispa}  to derive that
$$\overline{F}_{\varepsilon_n}(v_n)\geq \frac{1}{2\omega_n^2}\int_{\AA_{\varepsilon_n}\setminus \cup B_i^n}|j( v_n)|^2 - \zeta_{\max}\mu(\AA)+o(1)\,.$$
Setting $\tilde j_n(x)$ as in \eqref{deftildej} (with $\AA_{\varepsilon_n}$ in place of $\DD_{\varepsilon_n}$), up to a subsequence we have $\tilde j_n\rightharpoonup j_*$ weakly 
in $L^2(\AA;\RR^2)$ as $n\to+\infty$. By lower semicontinuity, we obtain
\begin{equation}\label{infhole1}
\liminf_{n\to+\infty} \overline{F}_{\varepsilon_n}(v_n)\geq \frac{1}{2}\int_{\AA}|j_*|^2 - \zeta_{\max}\mu(\AA)\,.
\end{equation}
In addition, arguing as in the proof of \eqref{curlj}, we deduce that
$$\curl\,j_*=\mu\quad\text{in $\mathscr{D}^\prime(\AA)$}\,,$$
and thus $\mu\in H^{-1}(\mathcal{D})$ since $j_*\in L^2(\AA;\RR^2)$. 

Next we introduce $h_\mu\in H_0^1(\mathcal{A})$ to be  the unique solution of 
$$\begin{cases} 
-\Delta h_\mu=\mu &\text{in $H^{-1}(\mathcal{A})$}\,,\\
h_\mu=0 & \text{on $\partial\mathcal{A}$}\,,
\end{cases}
$$
and  we set
$$\bar h_\mu:= h_\mu + \frac{1}{{\rm cap}(\BB)}\bigg(\int_{\partial \BB}j_*\cdot\tau + \frac{\partial h_\mu}{\partial \nu}\bigg)\xi\,.$$
By construction, we have $\curl(j_*+\nabla^\perp\bar h_\mu)=0$ in $H^{-1}(\AA)$ and 
$\int_{\partial \BB}( j_* + \nabla^\perp \bar h_\mu)\cdot\tau=0$. 
Hence there exists $g\in H^1(\AA)$ such that $j_*+\nabla^\perp \bar h_\mu=\nabla g$. Arguing as in the proof of Proposition~\ref{mindfix}, we derive
\begin{equation}\label{infhole2}
\int_\AA |j_*|^2= \int_\AA|\nabla \bar h_\mu|^2 +\int_\AA |\nabla g|^2\geq   \int_\AA|\nabla \bar h_\mu|^2\,.
\end{equation}
Then using \eqref{fctcap} and ${h_\mu}_{|\partial\AA}=0$, we obtain $\int_\AA \nabla h_\mu\cdot \nabla\xi =0$, 
so that  
\begin{equation}\label{infhole3}
\int_\AA|\nabla \bar h_\mu|^2= \int_\AA|\nabla h_\mu|^2 + \frac{1}{{\rm cap}(\BB)}\bigg(\int_{\partial \BB}j_*\cdot\tau + \frac{\partial h_\mu}{\partial \nu}\bigg)^2\geq \int_\AA|\nabla h_\mu|^2\,.
\end{equation}
Finally, using the Green representation of $h_\mu$ we obtain 
\begin{equation}\label{infhole4}
\frac{1}{2}\int_\AA|\nabla h_\mu|^2\,dx=\frac{1}{2}\iint_{\AA\times\AA}G(x,y)\,d\mu(x)d\mu(y)\,,
\end{equation}
and the  conclusion follows gathering \eqref{infhole1},  \eqref{infhole2},  \eqref{infhole3} and  \eqref{infhole4}.    
\prbox


\subsection{The upper bound inequality}

\noindent{\bf Proof of Theorem \ref{thmhole}, claim{\it  (iii)}.} We present here the proof in the case where the measure $\mu\in  H^{-1}(\AA)$ is absolutely continuous with respect to $\mathcal{H}^1\res\Sigma$, 
and more precisely for $\mu$ of the form \eqref{abscontmeas} with a nonvanishing density function $f$.  The general case follows by approximation as already pursued in Section \ref{secupbd}. 

For such a measure $\mu$ we first  proceed exactly as in Step 1 of the proof of Proposition~\ref{limsupabscont}, and we refer to it for the notation. For each integer $n$, we consider the function $f_n$ as defined in 
\eqref{deffn} and we set $\hat \mu_n:=f_n\mathcal{L}^2\res\AA$. Then $\hat \mu_n\rightharpoonup\mu$ weakly* as measures on $\AA$ as $n\to+\infty$. 
Next we consider the solution $h_n$ of 
$$\begin{cases}
-\Delta h_n=\omega_n f_n & \text{in $\AA$}\,,\\
h_n=0 & \text{on $\partial\AA$}\,.
\end{cases}$$
Arguing exactly as in Step 2 of the proof of Proposition \ref{limsupabscont}, we derive that 
$$\frac{1}{2}\int_\AA |\nabla h_n|^2\,dx\leq \pi D_n |\ln\varepsilon_n| + \omega_n^2\bar I(\mu)+o(\omega_n^2)\,. $$
Next we introduce 
$$\bar h_n:= h_n + \frac{\kappa_n}{{\rm cap}(\BB)} \xi\quad\text{with }\kappa_n:=\int_{\partial \BB}\frac{\partial h_n}{\partial\nu} - 2\pi\left[\frac{1}{2\pi}\int_{\partial \BB}\frac{\partial h_n}{\partial\nu}\right]\,, $$
where $[\cdot]$ denotes the integer part. Noticing that $\kappa_n=O(1)$ we deduce
\begin{multline}\label{energhbar}
\frac{1}{2}\int_\AA |\nabla \bar h_n|^2\,dx= \frac{1}{2}\int_\AA |\nabla h_n|^2\,dx + \frac{\kappa_n^2}{2{\rm cap}(\BB)^2} \int_\AA |\nabla \xi|^2\,dx\\
\leq \pi D_n |\ln\varepsilon_n| + \omega_n^2\bar I(\mu)+o(\omega_n^2)\,. 
\end{multline}
In view of \eqref{capa} we have 
$$\frac{1}{2\pi} \int_{\partial \BB}\frac{\partial \bar h_n}{\partial \nu}= \left[\frac{1}{2\pi}\int_{\partial \BB}\frac{\partial h_n}{\partial\nu}\right] \in \mathbb{Z}\,,$$
and since $\xi$ is harmonic in $\AA$, 
$$\frac{1}{2\pi} \int_{\partial B_{\varepsilon_n}(a^n_k)}\frac{\partial \bar h_n}{\partial \nu}= \frac{1}{2\pi} \int_{\partial B_{\varepsilon_n}(a^n_k)}\frac{\partial h_n}{\partial \nu}=-1\quad \text{for every $k=1,\ldots,D_n$}\,.$$
Hence, for any smooth Jordan curve $\Theta$ 
inside $\AA\setminus\cup_k \overline B_{\varepsilon_n}(a_k^n)$,
$$
\frac{1}{2\pi}\int_\Theta \nabla^\perp\bar h_n\cdot\tau \, \in 
  \mathbb{Z}\,,
$$
where $\tau:\Theta\to\mathbb{S}^1$ is any  smooth  vector field tangent to $\Theta$.  
Consequently, we may  locally define a phase $\phi_n$ in 
$\AA\setminus\cup_k \overline B_{\varepsilon_n}(a_k^n)$
by
$$   \nabla \phi_n(x) = -\nabla^\perp \bar h_n(x)\,,
    \quad x\in \AA\setminus\cup_k \overline B_{\varepsilon_n}(a_k^n)\,,$$
and then the map  
$\exp (i\phi_n(x))$ is well defined for every $x\in\AA\setminus\cup_k \overline B_{\varepsilon_n}(a_k^n)$. 

Finally we consider a profile function $\rho_n$ as defined in \eqref{defprof} and we set 
$$v_n(x):=\begin{cases}
\displaystyle \rho_n(x) e^{i\phi_n(x)} & \text{for $x\in \AA\setminus\cup_k \overline B_{\varepsilon_n}(a_k^n)$}\,,\\
0 & \text{otherwise}\,.
\end{cases}
$$
Since $-\Delta \bar h_n=\omega_n \hat \mu_n$ in $\AA$, using  \eqref{energhbar} we may proceed as in the proof of Proposition~\ref{limsupabscont} Step~3,  to prove that  
 $$\omega_n^{-1}{\rm curl}\,j(v_n)\to \mu \quad\text{strongly in $(C^{0,1}_0(\AA))^*$ as $n\to+\infty$}\,,$$
and that
\begin{equation}\label{upbdglhole}
E_{\varepsilon_n}(v_n,\AA)\leq \pi D_n |\ln\varepsilon_n| + \omega_n^2\bar I(\mu)+o(\omega_n^2)\,.
\end{equation}
To evaluate the rotation part of the energy, we first argue as for \eqref{U5} to obtain 
\begin{equation}\label{rothole1}
\Omega_n\int_\AA\nabla\zeta^\perp\cdot j(v_n)\,dx=-\pi D_n |\ln\varepsilon_n| -\zeta_{\rm max}\mu(\AA)\omega_n^2+o(\omega_n^2)\,. 
\end{equation}
In a similar way we derive that
\begin{equation}\label{rothole2}
\alpha_n\int_\AA \nabla^\perp \xi\cdot j(v_n)\,dx=-\alpha_n\int_\AA \nabla \bar h_n\cdot \nabla \xi\,dx+o(1)=-\alpha_n\kappa_n+o(1)=O(1)\,, 
\end{equation}
since $\alpha_n=O(1)$ and $\kappa_n=O(1)$. Then the conclusion follows gathering \eqref{upbdglhole}, \eqref{rothole1} and \eqref{rothole2}. 
\prbox

\subsection{Application to $J_\varepsilon$-$\,$minimizers}

\noindent{\bf Proof of Corollay \ref{asymptminhole}.} Let $\varepsilon_n\to 0^+$ be an 
arbitrary sequence. As for Corollary \ref{asymptmin}, the existence of minimizers for $J_{\varepsilon_n}$ and 
$\overline F_{\varepsilon_n}$ is classical, we omit the details. 
Then for every $n\in \NN$, let $u_n\in H^1(\AA;\CC)$ be a minimizer of $J_{\varepsilon_n}$. 
Since $J_{\varepsilon_n}(u_n)\leq \min \mathcal{H}_{\Omega_n}=O(\Omega_n^2)$, we can apply Proposition~\ref{energdecomp} to the sequence $\{u_n\}_{n\in\NN}$ to infer that 
\begin{equation}\label{quasiinf}
\min \mathcal{H}_{\Omega_n}\geq J_{\varepsilon_n}(u_n)= \min \mathcal{H}_{\Omega_n}+\omega_n^2\overline F_{\varepsilon_n}(v_n)+o(1)\,,
\end{equation}
with $v_n:=\bar u^\star_nu_n$. On the other hand, for any minimizer $\tilde v_n \in H^1(\AA;\CC)$ of $\overline F_{\varepsilon_n}$, we have  
$\overline F_{\varepsilon_n}(\tilde v_n)\leq \overline F_{\varepsilon_n}(1)=0$. Hence Proposition~\ref{energdecomp}  yields 
\begin{equation}\label{quasisup}
J_{\varepsilon_n}(u_n)\leq J_{\varepsilon_n}(u^\star_n\tilde v_n)=\min \mathcal{H}_{\Omega_n}+\omega_n^2\min\overline F_{\varepsilon_n}+o(1)\,. 
\end{equation}
Combining \eqref{quasiinf} with \eqref{quasisup}, we infer that  $\omega_n^2\overline F_{\varepsilon_n}(v_n)\leq o(1)$ and that $\{v_n\}$ is a sequence of quasi-minimizers for $\{\overline F_{\varepsilon_n}\}$, {\it i.e.}, 
\begin{equation}\label{quasimin}
\overline F_{\varepsilon_n}(v_n)= \min\overline F_{\varepsilon_n}+o(1)
\end{equation}
as $n\to +\infty$. By Theorem \ref{thmhole} there exists a subsequence $\{\varepsilon_{n_k}\}$ such that 
$$\omega_{n_k}^{-1} j(v_{n_k}) \to \mu_0 $$
strongly in $(C^{0,1}_0(\AA))^*$ for some nonnegative Radon measure $\mu_0\in H^{-1}(\AA)$ supprted by $\Sigma$. 
Using \eqref{quasimin} together with claims {\it (ii)} and {\it (iii)} in Theorem \ref{thmhole}, we deduce as in the proof of Corollary \ref{asymptmin} that 
$\mu_0$ minimizes $\mu\mapsto \bar I(\mu)-\zeta_{\rm max}\mu(\AA)$ over all nonnegative Radon measures supported by $\Sigma$. The minimizer is again unique and given by 
$$ \mu_0= \frac{\zeta_{\rm max}}{2 \bar I_*}\,\bar \mu_*\,,$$
whence the convergence of $\omega_n^{-1}j(v_n)$ along the full sequence (recall that $\bar I_*=\bar I(\mu_*)$ and that $\bar \mu_*$ is the minimizer of $\bar I$ over all probality measures 
supported by $\Sigma$, see \cite{ST}). In addition, 
$$\lim_{n\to+\infty}\overline F_{\varepsilon_n}(v_n)=\bar I(\mu_0)-\zeta_{\rm max}\mu_0(\AA) =-\frac{\zeta^2_{\rm max}}{4\bar I_*}\,,$$
which combined with \eqref{quasiinf} and \eqref{minengH} yields \eqref{minenergexphole}. 

Next it remains to prove that  
$$\omega_{n}^{-1}\, \curl\,j(u_{n})\mathop{\longrightarrow}\limits_{n\to+\infty} \mu_0$$
strongly in $(C^{0,1}_0(\AA))^*$. 
Indeed, given an arbitrary function $\varphi\in C^{0,1}_0(\AA)$, we have 
\begin{align*}
\frac{1}{\omega_{n}}\int_{\AA} j(u_{n})\cdot\nabla^\perp \varphi &=\frac{1}{\omega_{n}}\int_{\AA} j(v_{n})\cdot \nabla^\perp \varphi + \frac{1}{\omega_{n}} \int_{\AA} |v_n|^2j(u^\star_{n})\cdot \nabla^\perp \varphi \\
& = \langle \mu_0,\varphi\rangle +   \frac{1}{\omega_{n}} \int_{\AA} (|v_n|^2-1)j(u^\star_{n})\cdot \nabla^\perp \varphi\,+o(1)\,,
\end{align*}
where we have used \eqref{reljac} and the fact that ${\rm curl}\,j(u^\star_{n})=0$ in $\AA$. Arguing as in the proof of Proposition \ref{energdecomp} we estimate
$$\left| \frac{1}{\omega_{n}} \int_{\AA} (|v_n|^2-1)j(u^\star_{n})\cdot \nabla^\perp \varphi\right|\leq C\frac{\varepsilon_n}{\omega_n}\|j(u_n^\star)\|_{L^\infty(\AA)}\sqrt{E_{\varepsilon_n}(v_n)}
\leq C\frac{\varepsilon_n\Omega_n^2}{\omega_n} =o(1)\,,$$
where the constant $C$ only depends on $\varphi$, and the proof is complete. 
\prbox



\subsection*{Acknowledgements}
This work was initiated while V.M. was 
visiting the Department of Mathematics and Statistics at McMaster University: he thanks S.A., L.B. and the 
whole department for their warm hospitality. S.A. and L.B. are supported by NSERC (Canada) Discovery Grants.


\end{document}